\documentclass [12pt,oneside] {article}
\usepackage{latexsym,amsfonts,amssymb,amscd,amsmath,amsthm,mathrsfs,stmaryrd} 
\usepackage{graphicx, MnSymbol}
\usepackage{tikz}
\usetikzlibrary{fit}
\usepackage [all]{xy}
\SelectTips{cm}{12}

\theoremstyle{plain}
\newtheorem {Thm} {Theorem}[section]
\newtheorem* {Thm*} {Theorem}
\newtheorem* {Prop*} {Proposition}
\newtheorem {Lem}[Thm] {Lemma}
\newtheorem {Prop}[Thm] {Proposition}
\newtheorem {Cor}[Thm] {Corollary}

\theoremstyle {definition}
\newtheorem {Def}[Thm] {Definition}

\newtheorem {Exa}[Thm] {Example}
\newenvironment{Pf}[1]{{\noindent\sc Proof #1:}}{\qed\\}
\newcommand {\specialmap} [4] {{#1\negmedspace : #2 #3 #4}}

\newcommand {\map} [3] {\specialmap {#1} {#2}{\to} {#3}}
\newcommand {\longmap} [3] {\specialmap {#1} {#2}{\longrightarrow} {#3}}
\newcommand {\isomap} [3] {\specialmap {#1} {#2}{\overset {\cong{\phantom{.}}} 
            {\longrightarrow}} {#3}}

\newcommand {\Hom} {\operatorname {Hom}}
\newcommand {\id} {\operatorname{id}}

\newcommand {\at}[1] {\arrowvert_{#1}}

\newcommand {\ind} {\operatorname{ind}}
\renewcommand {\(} {\left(}
\renewcommand {\)} {\right)}

\def\acts{\hspace{.1cm}{
        \setlength{\unitlength}{.32mm}
        \linethickness{.09mm}
        \begin{picture}(8,8)(0,0)
        \qbezier(7,6)(4.5,8.3)(2,7)
        \qbezier(2,7)(-1.5,4)(2,1)
        \qbezier(2,1)(4.5,-.3)(7,2)
        \qbezier(7,6)(6.1,7.5)(6.8,9)
        \qbezier(7,6)(5,6.1)(4.2,4.4)
        \end{picture}\hspace{.1cm}
        }}

\newcommand {\sub} {\subseteq}

\newcommand {\CC} {\mathbb C}
\newcommand {\Cent} {{\on{\mD im}}}
\newcommand {\on}[1] {\operatorname{#1}}

\newcommand {\DCoh} {{\on{\mD^+_{\!coh}}}}
\newcommand {\DCon} {{\on{\mD^+_{\!constr}}}}

\newcommand {\Funead} {{\on{\mathcal F\!un_{r.e.}}}}

\newcommand {\GG}[1] {{\frac{1}{|#1|}}}

\renewcommand {\ind}[2] {\on{ind}\arrowvert_{#1}^{#2}}

\newcommand {\mA} {{\mathcal{A}}}
\newcommand {\mB} {{\mathcal{B}}}
\newcommand {\mC} {{\mathcal{C}}}
\newcommand {\mD} {{\mathcal{D}}}
\newcommand {\mF} {{\mathcal{F}}}
\newcommand {\mG} {{\mathcal{G}}}
\newcommand {\mK} {{\mathcal{K}}}

\newcommand {\ob} {{\operatorname{ob}}}

\newcommand {\Sn} {{S_n}}

\newcommand {\tensor}{\otimes}

\newcommand {\ZZ} {\mathbb Z}

\newcommand {\inv}{^{-1}}

\newcommand {\mI} {\mathcal I}
\newcommand {\mL} {\mathcal L}

\newcommand {\mV} {\mathcal V}
\newcommand {\mW} {\mathcal W}

\newcommand {\ul}[1]{\underline{#1}}

\newcommand {\RR} {{\mathbb R}}

\renewcommand{\theta}{\vartheta}
\newcommand {\ttr} {\mathbb Tr}
\renewcommand {\ul}[1] {\underline{#1}}

\renewcommand {\tilde}[1] {\widetilde{#1}}
\newcommand{\lwr} {\smallint} 

\title{Inner products of 2-representations}
\author{Nora Ganter\thanks{Ganter was supported by NSF-grant DMS-0504539 and by a
  Centenary Fellowship by the Faculty of Science at the University
  of Melbourne. Work on this project began while Ganter was at Colby
  College, Maine.
  }\\ The University of Melbourne}
\date {\today}
\begin{document}
\maketitle
\begin{abstract}
  We define and calculate inner products of 2-re\-pre\-sen\-ta\-tions.
  Along the way, we prove that the categorical trace $\ttr(-)$ of
  \cite[Sec.3]{Ganter:Kapranov} is multiplicative with respect to
  various notions of categorical tensor product, and we identify the
  center of the category of equivariant objects $\mV^G$ of
  \cite[Sec.4.2]{Ganter:Kapranov}. We discuss applications, ranging
  from Schur's result about the number of irreducible projective representations
  to a formula for the Hochschild cohomology of a global quotient orbifold.
\end{abstract}
\tableofcontents
\section{Introduction and statement of results}
The idea of a {\em representation of a group by functors in a
  category} goes back to Grothendieck \cite[p.196]{Grothendieck:Tohoku}, and  
categorified versions of representation theory 
have since been developed 
by a number of authors.\footnote{See for instance \cite{Elgueta}, \cite{Crane:Yetter},
\cite{Barrett:Mackaay}, \cite{Ostrik}, \cite{Freed},
\cite{Baez:Baratin:Freidel:Wise}, \cite{Deligne:action},
\cite{Seidel:Thomas}, \cite{Ben-Zvi} or \cite{Khovanov}.} 
(We recall some of the basic definitions in Section \ref{sec:2reps}.)
The corresponding character theory was introduced in
\cite{Bartlett} and \cite{Ganter:Kapranov}. Key ingredient
of this categorified character theory
is the notion of categorical trace: let $F$ be an endofunctor of a small
category. Then the {\em categorical trace} of $F$ is defined as the
set of natural transformations  
$$\ttr(F)= \on{\mathcal N\!at}(\id,F).$$ 
More generally, if $F$ is 
a 1-endomorphism 
in a 2-category, then $$\ttr(F)=2\Hom(1,F).$$ 
Let $\varrho$ be a linear 2-representation of a finite group $G$.
Then its {\em categorical character}
$X_\varrho$ sends $g\in G$ to the $k$-vector space
$$
  X_\varrho(g) = \ttr(\varrho(g)).
$$ 
The categorical character of $\varrho$ comes equipped with
conjugation isomorphisms
$$
  \isomap{\psi_{h}}{X_\varrho(g)}{X_\varrho(h\inv gh)},
$$
one for each pair of elements of $G$ (see
\cite[Prop.4.10]{Ganter:Kapranov}).  
As $g$ and $h$ vary, the $\psi_h$ are compatible, endowing
$X_\varrho$ with the structure of a {\em categorical class function} on $G$. Such
categorical class functions were first introduced by
Luztig \cite{Lusztig}, who named them {\em class sheaves}. 

For a commuting pair $(g,h)$,
the map $\psi_h$ is an automorphism of
$X_\varrho(g)$. Assume that all the $X_\varrho(g)$ are finite dimensional.
Then the {\em 2-character}
$\chi_\varrho$ of $\varrho$ is the function
$$
  \chi_\varrho(g,h) = tr(\psi_h),
$$
defined on pairs of commuting elements of $G$. It satisfies
$$
  \chi_\varrho(s\inv gs, s\inv hs) = \chi_\varrho(g,h).
$$
Such conjugation invariant 
functions on commuting pairs are called {\em 2-class functions}. 

Much of our motivation for writing \cite{Ganter:Kapranov} came from 
a character theory discovered by Hopkins, Kuhn and Ravenel in the
context of stable homotopy theory \cite{Hopkins:Kuhn:Ravenel}. 
Although of entirely different origin, these Hopkins-Kuhn-Ravenel
characters 
exhibit strikingly similar features to ours.
For instance, we proved in \cite[Cor.6.12]{Ganter:Kapranov} that 
the 2-character of an induced
2-representation is described by the same formula as the effect of 
transfer on Hopkins-Kuhn-Ravenel characters.

In this paper we will define inner products of
2-re\-pre\-sen\-ta\-tions and calculate them in
terms of categorical 
characters and 2-characters. On the level of 2-characters we obtain
\begin{equation}
  \label{eq:inner-product-of-characters}
  \langle \chi,\xi\rangle_G = \GG G\sum_{gh=hg}\chi(g,h)\cdot\xi(g,h)  
\end{equation}
(see Corollary \ref{2-inner-product-Cor}). This agrees with the formula for the
Strickland inner product of Hopkins-Kuhn-Ravenel characters 
in \cite[Prop.1.6]{Ganter:thesis}. 
Our definition will be
the categorical analogue of the inner product
$$
  \langle V, W\rangle_G := \dim_k\(V\tensor W\)^G
$$
on the representation ring $R(G)$ (compare \cite{Strickland}). This is 
closely related to the usual inner product, sending 
$(V,W)$ to 
$$
  \dim_k\(\Hom_G(W,V)\) = \langle V,W^*\rangle_G.
$$

The paper is organized as follows: in Section \ref{tensor-Sec},
we let $g$ and $h$ be endofunctors of linear categories $\mV$ and
$\mW$, respectively. 
Under a
finiteness condition we construct an isomorphism
$$
  \mu\negmedspace : \ttr(g)\tensor\ttr(h)\cong \ttr(g\boxtimes h)
$$
(Theorem \ref{tensor-Thm}).
Here $\boxtimes$ is the tensor product of linear categories
defined in \cite[Sec.1]{Ganter:Kapranov:symcat}.
Similar isomorphisms are found in the abelian case (using the Deligne tensor
product, \cite{Deligne}) and in the pre-triangulated case
(using the Bondal-Larsen-Lunts tensor 
product \cite{Bondal:Larsen:Lunts}).
We then discuss examples. The map $\mu$ exists in the more general
context of monoidal 2-categories. For some of those,
$\mu$ is an isomorphism.

In Section \ref{V^G-Sec} we let $\varrho$ 
be an action of $G$ on a linear category $\mV$, and consider the
category $\mV^G$ of equivariant objects of $\varrho$. We introduce the
{\em twisted group algebra} 
$$
  R_\varrho = \bigoplus_{g\in G}\ttr(g)
$$
acting on each object of $\mV^G$.
We contruct an isomorphism of $k$-algebras
\begin{equation}
  \label{eq:Dim-V^G}
  \Cent \(\mV^G\) \cong R_\varrho^G.  
\end{equation}
Here $$\Cent = \ttr(1)$$ is the center (of the category $\mV^G$).
As a corollary, we obtain 
$$
  \dim\(\Cent\(\mV^G\)\) = \GG G\sum_{gh=hg}\chi_\varrho(g,h).
$$
The category $\mV^G$ is characterized by a universal property, which
makes sense in the 2-categorical setup. All our proofs go through in
the context of 2-categories. 

Given two linear 2-representations $V$ and $W$ of $G$, we
can now define their inner product as
$$
  \langle V,W\rangle_G = \Cent\(V \boxtimes W\)^G,
$$
and conclude that, if $\langle V,W\rangle_G$ exists, its $k$-dimension is
calculated by 
\eqref{eq:inner-product-of-characters}. 
In Section \ref{applications-Sec}, we discuss applications of
\eqref{eq:Dim-V^G}. One 
special case is
a Theorem by Schur, counting the isomorphism classes of irreducible
projective representations with a given cocycle $c$.
Another special case is as follows:
let $X$ be a smooth projective variety over $\CC$, acted upon by
$G$. Let $\varrho$ be the resulting 2-representation of $G$ in
$\mV\!\on{ar_\CC}$ (see \cite[2.4.(c)]{Ganter:Kapranov}). In this context
\eqref{eq:Dim-V^G} becomes a result about orbifold Hochschild
cohomology, namely
\begin{equation}
  \label{eq:Hochschild-isom}
  \on{HH^\bullet}([X/G]) \cong \bigoplus_{[g],\alpha}
    \on{HH^{\bullet-codim(X^g_\alpha)}} \(X^g_\alpha,\det N^g\).  
\end{equation}
Here $X^g_\alpha$ are the connected components of the fixed point set $X^g$
of $g$, and $N^g$ is the normal bundle of $X^g$ in $X$; 
further, 
$$\on{HH^\bullet}(Y,\mF):=\on{Ext_{Y\times Y}}(\Delta_*\mathscr
  O_Y,\Delta_*\mathcal F).$$ 

I understand that an additive version of the isomorphism
\eqref{eq:Hochschild-isom}, with a very different proof, was 
known to Andrei Caldararu. 

The main advantage of our approach 
is multiplicativity: the right-hand side of
\eqref{eq:Hochschild-isom} is canoically 
an algebra, namely the twisted group algegra of the 2-representation
$\varrho$, and \eqref{eq:Hochschild-isom} is an algebra map.
Similar isomorphisms can be found in \cite{Caldararu:derived},
\cite{Dolgushev:Etingof}, \cite{Ginzburg:Kaledin} and
\cite{Baranovsky}, and 
it is a common theme that the proof 
of multiplicativity tends to be quite difficult. We note that our
proof of that fact 
is very simple (see page \pageref{page:multiplicativity}). 

Some interesting questions about multiplicativity
remain, for instance, whether it is
possible to give an explicit formula for the product on the
right-hand side, involving the double fixed point sets $X^g\cap X^h$.

Further applications of our results turn up in
\cite{Ganter:Kapranov:symcat}, where we define symmetric and
exterior powers of categories and calculate their effects on
characters of 2-representations. 
\subsubsection*{Acknowledgements} 
This paper is 
my account of results that were, to a large extent,
found in collaboration with Mikhail Kapranov. It is part of a series
of sequels to
\cite{Ganter:Kapranov}. 
Many crucial ideas 
were Mikhail's, and his comments on older drafts have been very
helpful. It is a pleasure to thank him for all this. 
I would also like to thank Alex Ghitza, Ezra
Getzler, and Craig Westerland for numerous helpful conversations on
the topic. 
\subsection{String diagram notation}
We will use string diagram notation for 2-morphisms.
For a nice, short introduction to string diagram
notation, we refer the reader to \cite{Caldararu:Willerton}, for a
more comprehensive account, we refer the reader to \cite{Bartlett}.
What follows is a mini summary of what we will use. Our conventions
follow those of Caldararu and Willerton. Bartlett's conventions are
slightly different.\footnote{In his diagrams 2-morphisms flow
  downwards, while ours flow upwards.}

  \begin{center}
    \begin{tikzpicture}[inner sep=2pt, scale=.7]
      \tikzset{
      short/.style={rectangle, draw=black},
      }
      \node at (-4,-1.5) [name=Y2] {$Y$};
      \node at (-6,-1.5) [name=phi2, scale=1.5] {$\Downarrow$};
      \node at (-5.9,-1.5) [name=phi2, anchor=west] {$\phi$};
      \node at (-8,-1.5) [name=X2] {$X$};
      \node at (-6,0.1) [anchor=south] {$f$};
      \node at (-6,-3.1) [anchor=north] {$g$};
      \node at (-4.8,-1) [name=a] {};
      \node at (-4.1,-.7) [name=c] {};
      \node at (-4.8,-2) [name=b] {};
      \node at (-4.1,-2.3) [name=d] {};
      \node at (-6,0) [name=e] {};
      \node at (-6,-3) [name=f] {};
      \draw [thick] (X2.north east) .. controls (-7,0.5) and (-5,0.5) .. (Y2.north
      west);  
      \draw [thick] (X2.south east) .. controls (-7,-3.5) and (-5,-3.5) .. (Y2.south
      west);  
      \draw [thick] (Y2.north west) -- (a);
      \draw [thick] (Y2.north west) -- (c);
      \draw [thick] (Y2.south west) -- (b);
      \draw [thick] (Y2.south west) -- (d);
      \node[draw=gray, fit=(e) (f) (X2) (Y2)] {};

      \node at (1,-3) [name=f, anchor=north] {$f$};
      \node at (1,-3) [name=b] {};
      \node at (1,0) [name=g, anchor=south] {$g$};
      \node at (1,0) [name=a] {};
      \node at (1,-1.5) [name=phi, short] {$\phi$};
      \node at (-0.2,-1.5) [name=Y] {$Y$};
      \node at (2.2,-1.5) [name=X] {$X$};

      \draw (g.south) -- (phi.north);
      \draw (phi.south) -- (f.north);
      \node[draw=gray, fit=(a) (b) (X) (Y)] {};
      \node at (-3,-4.5) [gray, anchor=north, text width=11cm]
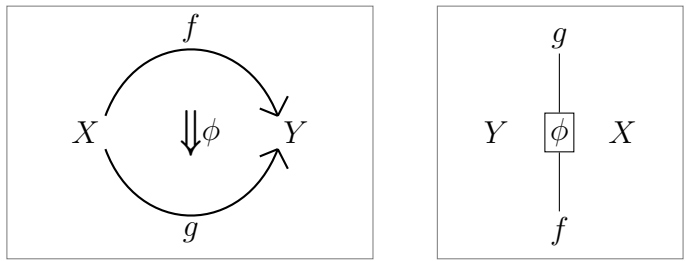
      {Traditional versus string diagram notation
        for a 2-morphism $\phi\negmedspace :f\Rightarrow
        g$ between 1-morphisms $f,g\negmedspace : X\to Y$.};
      \node at (0,-3) {};
    \end{tikzpicture}
  \end{center}
String diagrams are Poincare dual to the classical notation. They are
read from the right to the left (direction of 1-morphisms) and from
the bottom to the top (direction of 2-morphisms). 
The objects $X$ and $Y$ are sometimes omitted from the picture,
when they are clear from the context.

Identity morphisms are denoted by omission:
  \begin{center}
    \begin{tikzpicture}[inner sep=2pt, scale=.7]
      \tikzset{
      short/.style={rectangle, draw=black},
      }
      \node at (1,0) [name=f, anchor=south] {$f$};
      \node at (1,-1.5) [name=phi, short] {$\phi$};
      \node at (1,-3) [name=X] {$X$};
      \node at (0,0) [name=a] {};
      \node at (2,-3) [name=b] {};
      \node at (1,-4.3) [gray, anchor=north, text width=3cm] {A
        2-morphism $\phi\negmedspace :\id_X\Rightarrow f$};

      \draw (f.south) -- (phi.north);
      \node[draw=gray, fit= (b) (a)] {};

      \node at (8,-3) [name=f, anchor=north] {$f$};
      \node at (8,-1.5) [name=phi, short] {$\phi$};
      \node at (8,0) [name=X] {$X$};
      \node at (7,0) [name=a] {};
      \node at (9,-3) [name=b] {};
      \node at (8,-4.3) [gray, anchor=north, text width=3cm] {A
        2-morphism $\phi\negmedspace :f\Rightarrow \id_X$};

      \draw (f.north) -- (phi.south);
      \node[draw=gray, fit= (b) (a)] {};

      \node at (15,0) [name=f1, anchor=south] {};
      \node at (15,-3) [name=f2, anchor=north] {$f$};
      \node at (14,0) [name=a] {};
      \node at (16,-3) [name=b] {};
      \node at (15,-4.3) [gray, anchor=north, text width=3cm] {The identity\\
        2-morphism $\id_f$};

      \draw (f1.south) -- (f2.north);
      \node[draw=gray, fit= (b) (a)] {};
    \end{tikzpicture}
  \end{center}

Composition is denoted by juxtaposition:
\medskip
  \begin{center}
    \begin{tikzpicture}[inner sep=2pt, scale=.7]
      \tikzset{
      short/.style={rectangle, draw=black},
      }
      \node at (1.4,0) [name=h, anchor=south] {$h$};
      \node at (1.4,-1.5) [name=psi, short] {$\psi$};
      \node at (1.4,-3.5) [name=phi, short] {$\phi$};
      \node at (1.4,-5) [name=f, anchor=north] {$f$};
      \node at (1.4,-2.5) [name=g, anchor=west] {$g$};
      \node at (3.1,-2.5) [name=X] {$X$};
      \node at (-0.3,-2.5) [name=Y] {$Y$};
      \node at (-.1,-.5) [name=a] {};
      \node at (2.9,-4.5) [name=b] {};
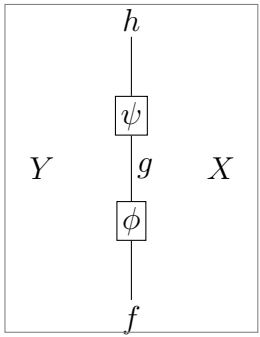
      \node at (1.4,-6) [gray, anchor=north, text width=3.5cm] {Vertical
      composition $\psi\circ_1\phi$};

      \draw (h.south) -- (psi.north);
      \draw (psi.south) -- (phi.north);
      \draw (phi.south) -- (f.north);
      \node[draw=gray, fit= (b) (a)] {};
 
      \node at (7,0) [name=k, anchor=south] {$k$};
      \node at (7,-1.5) [name=psi, short] {$\psi$};
      \node at (7,-5) [name=h, anchor=north] {$h$};
 
      \node at (9,0) [name=g, anchor=south] {$g$};
      \node at (9,-3.5) [name=phi, short] {$\phi$};
      \node at (9,-5) [name=f, anchor=north] {$f$};

      \node at (10.3,-2.5) [name=X] {$X$};
      \node at (5.7,-2.5) [name=Z] {$Z$};
      \node at (8,-2.5) {$Y$};
      \node at (6.3,-.5) [name=a] {};
      \node at (9.7,-4.5) [name=b] {};

      \draw (k.south) -- (psi.north);
      \draw (psi.south) -- (h.north);
      \draw (g.south) -- (phi.north);
      \draw (phi.south) -- (f.north);
      \node[draw=gray, fit= (b) (a)] {};

      \node at (13,0) [name=k, anchor=south] {$k$};
      \node at (13,-3.5) [name=psi, short] {$\psi$};
      \node at (13,-5) [name=h, anchor=north] {$h$};
 
      \node at (15,0) [name=g, anchor=south] {$g$};
      \node at (15,-1.5) [name=phi, short] {$\phi$};
      \node at (15,-5) [name=f, anchor=north] {$f$};

      \node at (16.3,-2.5) [name=X] {$X$};
      \node at (11.7,-2.5) [name=Z] {$Z$};
      \node at (14,-2.5) {$Y$};
      \node at (12.3,-.5) [name=a] {};
      \node at (15.7,-4.5) [name=b] {};

      \draw (k.south) -- (psi.north);
      \draw (psi.south) -- (h.north);
      \draw (g.south) -- (phi.north);
      \draw (phi.south) -- (f.north);
      \node[draw=gray, fit= (b) (a)] {};

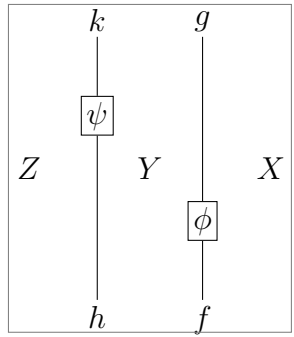
      \node at (11,-6) [gray, anchor=north, text width=8cm] {Two ways
        of drawing the horizontal composition of 2-morphisms,
        $\psi\phi\negmedspace : hf\Rightarrow kg$. The fact that both
        diagrams describe the same 2-morphisms is known as the ``butterfly
        identity'' $(\psi g)\circ_1(h\phi)
        =(k\phi)\circ_1(\psi f)=:\psi\phi$.};
    \end{tikzpicture}
  \end{center}

\begin{Def}
  An {\em adjunction} in a 2-category $\mC$ consists of the following data:
  two objects, $X$ and $Y$, 1-morphisms $\map fXY$ and $\map gYX$ and
  2-morphisms 

  
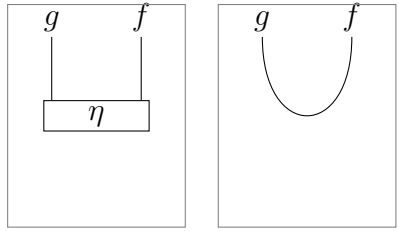
\begin{center}
    \begin{tikzpicture}[inner sep=2pt, scale=.7]
      \tikzset{
      long/.style={rectangle, draw=black, minimum width=1.4cm}
      }
      \node at (.65,0) [name=g, anchor=south] {$g$};
      \node at (2.35,0) [name=f, anchor=south] {$f$};
      \node at (1.5,-1.5) [name=eta, long] {$\eta$};
      \node at (.5,-.2) [name=a] {};
      \node at (2.5,-2.8) [name=b] {};

      \draw (g.south) -- ([xshift=-.85cm]eta.north);
      \draw (f.south) -- ([xshift=.85cm]eta.north);
      \node[draw=gray, fit= (b) (a)] {};

      \node at (4.65,0) [name=g, anchor=south] {$g$};
      \node at (6.35,0) [name=f, anchor=south] {$f$};
       \node at (4.5,-.2) [name=a] {};
      \node at (6.5,-2.8) [name=b] {};

      \draw (g.south) .. controls (4.65,-2) and (6.35,-2) .. (f.south);
      \node[draw=gray, fit= (b) (a)] {};  

      \node at (3.5,-4) [gray, anchor=north, text width=5cm]
      {$\eta\negmedspace :\id_X\Rightarrow gf$, called the {\em unit}
        and denoted by either of
        these two pictures};

      \node at (9.65,-3) [name=f, anchor=north] {$f$};
      \node at (11.35,-3) [name=g, anchor=north] {$g$};
      \node at (10.5,-1.5) [name=epsilon, long] {$\epsilon$};
      \node at (9.5,-.2) [name=a] {};
      \node at (11.5,-2.8) [name=b] {};

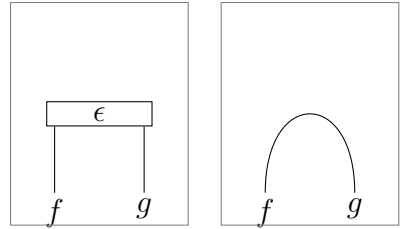
      \draw (g.north) -- ([xshift=.85cm]epsilon.south);
      \draw (f.north) -- ([xshift=-.85cm]epsilon.south);
      \node[draw=gray, fit= (b) (a)] {};

      \node at (13.65,-3) [name=f, anchor=north] {$f$};
      \node at (15.35,-3) [name=g, anchor=north] {$g$};
      \node at (13.5,-.2) [name=a] {};
      \node at (15.5,-2.8) [name=b] {};

      \draw (f.north) .. controls (13.65,-1) and (15.35,-1) .. (g.north);
      \node[draw=gray, fit= (b) (a)] {};  
      
      \node at (12.5,-4) [gray, anchor=north, text width=5cm] {and
      $\epsilon\negmedspace : fg\Rightarrow \id_Y$, called the {\em
        counit} and denoted by either
      of these two pictures.};
    \end{tikzpicture}
  \end{center}
  These are required to satisfy
  
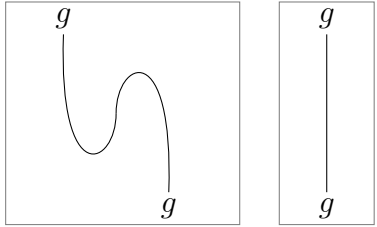
\begin{center}
    \begin{tikzpicture}[inner sep=2pt, scale=.7]
      \tikzset{
      long/.style={rectangle, draw=black, minimum width=1.4cm}
      }
      \node at (0,0) [name=g1, anchor=south] {$g$};
      \node at (2,-3) [name=g2, anchor=north] {$g$};
      \node at (1,-1.5) [name=f] {};
      \node at (-.25,-.2) [name=a] {};
      \node at (2.5,-2.8) [name=b] {};

      \draw (g1.south) .. controls (-.1,-3) and (1,-2.5)
      .. (f.center);
      \draw (f.center) .. controls (1,-.5) and (2.1,0) .. (g2.north);
      \node[draw=gray, fit= (b) (a)] {};

      \node at (5,0) [name=g1, anchor=south] {$g$};
      \node at (5,-3) [name=g2, anchor=north] {$g$};

       \node at (4.55,-.2) [name=a] {};
      \node at (5.45,-2.8) [name=b] {};

      \draw (g1.south) -- (g2.north);
      \node[draw=gray, fit= (b) (a)] {};  

      \node at (2.6,-4) [gray, anchor=north, text width=5cm]
      {$(\id_g\epsilon)(\eta\id_g)=\id_g$, in other words, the
        2-morphisms in these two
        string diagrams are required to be equal,};

      \node at (12,0) [name=f1, anchor=south] {$f$};
      \node at (10,-3) [name=f2, anchor=north] {$f$};
      \node at (11,-1.5) [name=g] {};
      \node at (9.75,-.2) [name=a] {};
      \node at (12.5,-2.8) [name=b] {};

      \draw (f1.south) .. controls (12.1,-3) and (11,-2.5)
      .. (g.center);
      \draw (g.center) .. controls (11,-.5) and (9.9,0) .. (f2.north);
      \node[draw=gray, fit= (b) (a)] {};

      \node at (15,0) [name=g1, anchor=south] {$f$};
      \node at (15,-3) [name=g2, anchor=north] {$f$};

       \node at (14.55,-.2) [name=a] {};
      \node at (15.45,-2.8) [name=b] {};

      \draw (g1.south) -- (g2.north);
      \node[draw=gray, fit= (b) (a)] {};  

      \node at (12.6,-4) [gray, anchor=north, text width=5cm]
      {and $(\epsilon\id_f)(\id_f\eta)=\id_f$,\\ i.e., the
        2-morphisms in these two
        string diagrams are required to be equal.};
    \end{tikzpicture}
  \end{center}
  We say ``$f$ is left adjoint to $g$''.
\end{Def}
For more on adjunctions in 2-categories, see for instance \cite[p.\ 13
    ff.]{Caldararu:Willerton}. 
\section{Multiplicativity of the categorical trace}
\label{tensor-Sec}   
Let $k$ be an algebraically closed field. By a {\em linear category}, we
will always mean an additive $k$-linear category, and by a {\em linear
  functor} we mean one that preserves this structure. 
In \cite{Ganter:Kapranov:symcat}, we consider three different notions of
tensor product of categories:
\begin{enumerate}
\item if $\mV$ and $\mW$ are linear categories, we consider the
  uncompleted tensor product 
  $\mV\boxtimes \mathcal W,$ 
  following \cite{Kapranov:Voevodsky} and \cite{Bakalov:Kirillov},
\item if $\mathcal A$ and $\mathcal B$ are abelian linear categories,
  we consider the completed 
  version $\mA\widehat\boxtimes \mathcal B$, following Deligne
  \cite{Deligne},
\item if $\mC$ and $\mathcal D$ are perfect dg-categories over $k$, we consider
  the
  completed tensor pro\-duct $\mC\widehat \boxtimes \mathcal D$,
  following \cite{Bondal:Larsen:Lunts}.
\end{enumerate}
The goal of this section is to find isomorphisms
$$
  \mu\negmedspace : \ttr(g)\tensor\ttr(h)\cong \ttr(g\boxtimes h)
  \cong \ttr(g\widehat \boxtimes h),
$$
where $g$ and $h$ are endofunctors of the relevant categories. We can
always define a map $\mu$, but in order to show that it is an isomorphism,
we need a finiteness condition.
\subsection{The multiplicativity theorem}
The first two constructions above are characterized by similar universal
properties. For instance, the data of (1) 
consist of a linear category $\mV\boxtimes\mW$
together with a $k$-bilinear functor
$$
  \longmap\boxtimes{\mV\times\mW}{\mV\boxtimes \mW},
$$
which is universal in the following sense:
let $\mathcal Z$ be a third linear category.
Then precomposition with $\boxtimes$ defines an equivalence of categories
$$
  \mF\!\on{un}_{\on{lin}}(\mV\boxtimes \mW,\mathcal Z)
  \stackrel\simeq\longrightarrow
  \mF\!\on{un}_{\on{bil}}(\mV\times \mW,\mathcal Z).
$$
Here $\mF\!\on{un}_{\on{lin}}$ and
$\mF\!\on{un}_{\on{bil}}$ stand for the categories of $k$-linear and
$k$-bilinear functors.   

The categorical tensor products in (1) and (3) always exist. In the
abelian case, (2),
Deligne has proved the existence of $\mA\widehat\boxtimes\mB$ under the
following finiteness 
condition: 
\begin{Def}[{compare \cite[(2.12.1)]{Deligne}}] 
  We say that a linear abelian category $\mA$ {\em satisfies Deligne's
    finiteness 
    condition}, if all objects of $\mA$ have finite length and all
  $\Hom$-sets are finite dimensional $k$-vector spaces.
\end{Def}
Kapranov and Voevodsky argued that $\boxtimes$ endows the 2-category
$\mC\!\on{at_k}$ of $k$-linear categories with the structure of a  monoidal
2-category.
The reader interested in the full twelve pages of axioms
is referred to
\cite{Kapranov:Voevodsky}.\footnote{The definition was later revised
  by Crans \cite{Crans}, see also \cite{Baez:Langford}.}
The following proposition, which follows immediately from the
universal propetry, summarizes what we will need.
\begin{Prop}
(1)
  Let $\mV$, $\mV'$, $\mW$ and $\mW'$ be linear categories. Then there
  exist a bilinear functor
  \begin{eqnarray*}
  (-\boxtimes-)\negmedspace :  \mF\!\on{un_{lin}}(\mV,\mV')\times \mF\!\on{un_{lin}}(\mW,\mW')
    &\longrightarrow
    &\mF\!\on{un_{lin}}(\mV\boxtimes\mW,\mV'\boxtimes\mW')\\
  \end{eqnarray*}
  and a natural isomorphism
  $$
    \iota\negmedspace : \boxtimes\circ(-,-)\Longrightarrow 
    (-\boxtimes -)\circ\boxtimes.
  $$
  The pair $\((-\boxtimes-),\iota\)$ is unique up to unique natural
  isomorphism (preserving $\iota$).
  We have canonical functor isomorphisms
  $$
    \isomap{\kappa_{g,f,h,k}}{(gf)\boxtimes (hk)}{(g\boxtimes h)
      (f\boxtimes k)}, 
  $$
    whenever the compositions on the left-hand side exist. They are natural
    in $g$, $f$, $h$, and $k$.
    Further, there is a functor isomorphism
  $$
    \isomap{\kappa_1}{\id_\mV\boxtimes\id_\mW}{\id_{\mV\boxtimes\mW}}.
  $$
  These $\kappa_?$ are compatible with $\iota$ in an obvious sense and
  make the usual associativity pentagon and unit diagrams commute.

\medskip
(2) Let $(\mV\widetilde\boxtimes\mW,\widetilde\boxtimes)$ be a second
pair satisfying the universal property of
$(\mV\boxtimes\mW,\boxtimes)$. Then there exist an equivalence  
$$
  \longmap E{\mV\boxtimes\mW}{\mV\widetilde\boxtimes\mW}
$$
and a natural isomorphism
$$
  \iota_E\negmedspace : \widetilde\boxtimes \Rightarrow
  E\circ\boxtimes. 
$$ 
The pair $(E,\iota_E)$, is unique up to unique natural
isomorphism
(preserving $\iota_E$).
\end{Prop}
%

A similar statement holds for the abelian case and right-exact linear
functors.

\begin{Def}\label{def:mu}
  Let $\mV$ and $\mW$ be linear categories, let $g$ be an endofunctor
  of $\mV$, and let $h$ be an endofunctor of $\mW$. Then the proposition
  gives a bilinear map
  $$
    \longmap{\widetilde\mu}{\ttr(g)\times
      \ttr(h)} {    \on{\mathcal
        N\!at}(\id_\mV\boxtimes\id_\mW,g\boxtimes h)}\longrightarrow
    \ttr(g\boxtimes h)
  $$
  $$(\phi,\psi)\longmapsto\(\phi\boxtimes\psi\)\circ\kappa_1\inv.$$
  We let $\mu$ be the resulting linear map
  $$
    \longmap\mu{\ttr(g)\tensor\ttr(h)}{\ttr(g\boxtimes h)}.
  $$
\end{Def}
Note that $\mu$ is natural in $g$ and $h$.
The following is a reformulation of
\cite[Prop.5.13 (i),(v)]{Deligne}: 
\begin{Prop}\label{completion-map-Prop}
  Let $\mA$ and $\mB$ be $k$-linear abelian categories satisfying
  Deligne's finiteness condition. Then the completed tensor product
  $\mA\widehat\boxtimes\mB$ exists. Further, there is a 
  fully faithful ``completion'' functor
  $$
    \mA\boxtimes\mB\longrightarrow  \mA\widehat\boxtimes\mB,
  $$
  obtained by applying the universal property of $\boxtimes$ to
  $\widehat\boxtimes$. 
\end{Prop}
  In particular, the essential
  image of
  $\map{\widehat\boxtimes}{\mA\times\mB}{\mA\widehat\boxtimes\mB}$ is
  contained in $\mA\boxtimes\mB$, and 
  the two functors $\widehat\boxtimes$ and $\boxtimes$ are identified
  when viewed as functors onto their essential images. 
  
  By construction, there is a similar fully faithful ``completion
  functor'' for the Bondal-Larsen-Lunts tensor product 
  of pre-triangulated categories. 

\begin{Thm}
\label{tensor-Thm}
  (1) Let $\mV$ and $\mW$ be linear categories
  such that all $\Hom$-sets are finite dimensional $k$-vector spaces.
  Let $g$ be a linear en\-do\-func\-tor of $\mV$, and
  let $h$ be a linear en\-do\-func\-tor of $\mathcal W$. 
  Assume that either $\ttr(g)$ or $\ttr(h)$ is finite dimensional.
  Then the map $\mu$ of Definition \ref{def:mu} is an isomorphism.
  
  \medskip
  (2) Let $\mA$ and $\mathcal B$ be linear abelian categories satisfying
  Deligne's finiteness condition.
   Let $g$ be a right-exact endofunctor of $\mA$, and let $h$ be a
  right-exact endofunctor of $\mathcal B$. 
  Assume that either $\ttr(g)$ or $\ttr(h)$ is finite dimensional.
  Then we have an isomorphism
  $$
    \ttr(g\boxtimes h)\longrightarrow
    \ttr(g\widehat \boxtimes h),
  $$
  which is natural in $g$ and $h$.
\end{Thm}
\begin{Pf}{}
  Part (2):
  let $\mA$, $\mathcal B$, $g$ and $h$ be as in (2).
  Then we have 
  \begin{eqnarray*}
    \ttr (g\widehat\boxtimes h) 
    &\cong& \on{\mathcal N\!
      at}(\widehat\boxtimes,\widehat\boxtimes\circ(g\times h))\\
    &=& \on{\mathcal N\!
      at}(\boxtimes,\boxtimes\circ(g\times h))\\
    &\cong& 
     \ttr(g\boxtimes h).
  \end{eqnarray*}   
  Here the first and third equality come from the universal
  properties of $\widehat \boxtimes$ and $\boxtimes$. The key step is
  the second equality, which follows from Proposition 
  \ref{completion-map-Prop}. 

  Part (1):
  Let $(X,Y)$ be an object of $\mV\times\mathcal W$.
  By construction of $\mV\boxtimes\mW$, we have
  $$
    \Hom_{\mathcal V\boxtimes \mathcal W}(X\boxtimes
    Y,gX\boxtimes hY)\cong
    \Hom_{\mathcal V}(X,gX)\tensor_k  \Hom_{\mathcal W}(Y,hY)
  $$
  (see \cite[Sec.1]{Ganter:Kapranov:symcat}).
  Under this identification,
  $$
    \mu\(\sum_i\phi_i\tensor\psi_i\)_{(X,Y)} =\sum_i
    \phi_{i,X}\tensor\psi_{i,Y}. 
  $$
  To define an inverse of $\mu$,  
  let $$\eta\negmedspace :\boxtimes\Rightarrow
  \boxtimes\circ(g,h)$$ be a natural transformation. Fix an object $Y$
  of $\mathcal W$, and choose a 
  basis $B^Y=(\beta^Y_1,\dots,\beta_m^Y)$ of $\Hom_{\mathcal W}(Y,hY)$.
  For each $X\in\ob\mV$, there is a unique way to express
  $\eta_{(X,Y)}$ in the form
  $$
    \eta_{(X,Y)} =: \sum_{i=1}^m\phi^X_i\tensor \beta^Y_{i}.
  $$
  As $X$ varies, naturality of $\eta$ implies that, for each $i$, the
  maps $\phi_i^X$ form a natural transformation $\phi_i$ from $\id_\mV$ to $g$.
  So, 
  $$
    \eta_Y := \sum_{i=1}^m\phi_i\tensor \beta_i^Y
  $$
  is a well-defined element of
  $$\ttr_\mV(g)\tensor\Hom_{\mW}(Y,hY).$$
  Without loss of generality, we assume $\ttr(g)$ to be finite dimensional.
  Then we may pick a $k$-basis of $\ttr_\mV(g)$ and repeat our
  argument. This yields an element 
  $$
    \nu(\eta)\in\ttr_\mV(g)\tensor_k\ttr_{\mW}(h).
  $$
  By construction, $\nu$ is a right-inverse to $\mu$.
  To see that $\nu$ is also a left-inverse, we pick $\phi\in\ttr(g)$
  and $\psi\in\ttr(h)$ and calculate
  $\nu(\mu(\phi\tensor\psi))$. Fix $Y$, and write $\psi_Y$ as a
  linear combination of the basis $B^Y$,
  $$\psi_Y =:
  \sum_ia_i\beta_i^Y.$$  
  In the construction of $\nu$, this implies
  $\phi_i^X = a_i\phi_X$, hence $\phi_i = a_i\phi$ and  
  \begin{eqnarray*}
    \eta_Y &=& \sum_i \phi_i\tensor \beta_i \\ 
    &=& \sum_i a_i\phi\tensor \beta_i \\
    &=& \phi\tensor \sum_i a_i \beta_i \\
    &=& \phi\tensor \psi_Y.   
  \end{eqnarray*}
  Hence $$\nu(\mu(\phi\tensor\psi))=\phi\tensor\psi,$$
  as claimed, and we have proved that $\nu$ is a left-inverse of $\mu$.
\end{Pf}{}
\subsection{Examples}
\nopagebreak[4]
\subsubsection{2-Vector spaces}
Let $\mV$ be a semisimple 
linear category 
with finitely many isomorphism classes of simple objects. 
Let $m$ be the number of isomorphism classes of
simple objects in $\mV$. Then there is an equivalence of categories
$$
  \mV\simeq \on{Vect}_k^m.
$$
This is the reason why categories like $\mV$ are called {\em $2$-vector
  spaces} \cite{Kapranov:Voevodsky}.
If $\mV$ and $\mW$ are 2-vector spaces, Osorno has independently found
an isomorphism 
$$\ttr(g\boxtimes h) \cong \ttr(g)\tensor\ttr(h)$$ \cite{Osorno}. 
We give a brief summary of her argument:
a linear functor between two 2-vector spaces can be represented by
a matrix with entries in $\on{Vect}_k$. Its categorical trace
is the direct sum of the diagonal entries. 
Given equivalences $\mV \simeq \on{Vect}_k^m$ and $\mW \simeq
\on{Vect}_k^n$, one obtains an equivalence
$$
  \mV\boxtimes \mW\simeq \on{Vect}_k^{mn}.
$$
Let $g$ be a linear endofunctor of $\mV$, and let $h$ be a linear
endofunctor of $\mW$.
As one would expect, the matrix for $g\boxtimes h$ is the ``Kronecker
product'' of the matrices for $g$ and for $h$. 
The classical computation for 
$$
  tr(g)\cdot tr(h) = tr(g\tensor h)
$$
goes through and shows that 
$\ttr(g\boxtimes h)$ is isomorphic to 
$\ttr(g)\tensor\ttr(h)$. 
\subsubsection{Algebras}\label{tensor-algebra-Sec}
Let $A$ be an associative and unitary finite dimensional 
$k$-algebra,
%
%
and let $\on{A-mod^f}$ 
be the category of finitely presented right
$A$-modules.\footnote{This is the category denoted $(A)_{\on{coh}}$ 
in \cite{Deligne}.}
Let $g$ be a right-exact endofunctor of $\on{A-mod^f}$. 
Recall (e.g.\ from
\cite{Deligne}) that $g$ is uniquely determined by the
$A$-$A$-bimodule $M:= g(A)$. Here the left module structure on $M$ comes
from the left action of $A$ on itself:
$$
  A\longrightarrow \on{End}(A)\longrightarrow \on{End}(M).
$$
The categorical trace of $g$ is isomorphic to the center of $M$:
  \begin{eqnarray*}
    \ttr(g) 
    &\cong&\Hom_{\on{A-A}}(A,g(A))\\
    & =& \on{Center}_A(M).
  \end{eqnarray*}
Here $\Hom_{\on{A-A}}$ stands for morphisms of $A$-$A$-bimodules, and
the second isomorphism identifies $f\in \Hom_{\on{A-A}}(A,M)$ with the element
$f(1)$ in $M$. 

Let $B$ be a second finite dimensional $k$-algebra.
Let $h$ be a right-exact
linear endofunctor of $\on{B-mod^f}$.
Recall from \cite[Prop.5.3]{Deligne} that we have an equivalence of
abelian categories 
\begin{equation}
  \label{A-Coh-tensor-Eqn}
  \(\on{A-mod^f}\)\widehat\boxtimes\(\on{B-mod^f}\)\simeq
  \on{(A\tensor B)-mod^f},   
\end{equation}
where $\widehat\boxtimes$ is the Deligne tensor product. 

  Let $M:=g(A)$ and $N:= h(B)$.
  Under the equivalence \eqref{A-Coh-tensor-Eqn}, the isomorphism
  $$
    \longmap{\mu}{\ttr(g)\tensor\ttr(h)}{\ttr(g\widehat\boxtimes h)}
  $$
  of Theorem \ref{tensor-Thm} is
  identified with the canonical map
  $$
    \longmap{\mu}{\on{Center_A}(M)\tensor \on{Center_B}(N)}  
    {\on{Center_{A\tensor B}}(M\tensor N)}.
  $$  
  The fact that $\mu$ is an isomorphism is the degree zero part of
  the K\"unneth theorem for Hochschild
  cohomology (c.f.\ \cite[X.7.4]{MacLane}).

There is, of course, overlap between this example and the previous one:
  if $A$ is semisimple, then $\on{A-mod^f}$ is a 2-vector space. So,
  $$\on{A-mod^f}\simeq \on{Vect_k}^r,$$ where $r$ is the number
  of isomorphism classes of simple $A$-modules. 
\begin{Exa}
  Take $g=\id$ to be the
  identity functor. We have seen two different ways to calculate the
  categorical trace of $\id$: on one hand,
  $$\ttr(\id)\cong\on{Center}(A).$$ 
  On the other hand, we may
  calculate $\ttr(\id)$ as the direct sum of the diagonal entries of
  the $r\times r$ identity 2-matrix, $$\ttr(\id)\cong k^{\oplus r}.$$
  So, we have recovered
  the well known fact that $$\dim(\on{Center}(A)) = r.$$    
\end{Exa}

More generally, let $M$ be an $A$-$A$-bimodule. Then the same argument yields
  $$\dim(\on{Center_A}(M)) = \sum m_i,$$ where
  $m_i$ is the number of times that the simple object
  $W_i$ occurs as summand inside $W_i\tensor 
  M$. Here the sum runs over all isomorphism classes of simple objects
  of $\on{A-mod^f}$. 
\subsubsection{Coherent sheaves}
Let $X$ be a smooth projective variety over $k$, and let $\mC\!\on{oh}(X)$
be the category of coherent sheaves on $X$.
Let $\DCoh(X)$ be the derived category of bounded below complexes in
$\mC\!\on{oh}(X)$. This is a triangulated category, and it possesses
an enhancement $\mI(X)$ 
such that, for any two varieties $X$ and $Y$ as above, the completed
tensor product  
$$\mI(X)\medspace \widehat\boxtimes \medspace\mI(Y)$$ is quasi-equivalent to
$\mI(X\times Y)$ (see \cite[Th. 5.5]{Bondal:Larsen:Lunts} compare also
\cite[(1.7.4)(c)]{Ganter:Kapranov:symcat}). 
Let $g$ be an automorphism of $X$, and let $h$ be an automorphism of
$Y$. Then we have the direct image functors $g_*$ and $h_*$ 
acting on, respectively, $\mI(X)$ and $\mI(Y)$. The above equivalence
identifies $g_*\widehat\boxtimes h_*$ with $(g\times h)_*$.
Applying Theorem \ref{tensor-Thm}, we conclude
$$
  \ttr^\bullet(g_*)\tensor\ttr^\bullet(h_*) \cong
  \ttr^\bullet((g\times h)_*).
$$
These categorical traces in $\mI(X)$ are hard to get a handle
on. Instead, one often works with
bicategories, compare to Section \ref{sec:DCoh-mu}. 
\section{Traces in monoidal 2-categories}
By a {\em monoidal 2-category} we will mean a $k$-linear semistrict monoidal
2-category in the sense of \cite{HDAI}.
In this setup, the axioms still imply the existence of $g\boxtimes h$,
unique up to a specified isomorphism and the existence of a map
$$\map\mu{\ttr(g)\tensor\ttr(h)}{\ttr(g\boxtimes h)}.$$
In general, there is no reason for $\mu$ to be an isomorphism, but
in the following we will see some examples where it is. 
To be precise, the lax 2-categories in the following sections would need to be
strictified in an appropriate sense (see \cite[2.12,
4.3]{Kapranov:Voevodsky}) in order for the formalism of \cite{HDAI} to
apply. Rather than working out the precise formalism in the lax case,
we choose an ad hoc approach to the individual examples.
\subsection{Bimodules, strict version}
\label{sec:bimodules,strict}
We recall 
from \cite[2.2(c)]{Ganter:Kapranov} the k-linear 2-category $\mathcal
B\!\on{im}_k$.
The objects of $\mathcal B\!\on{im}_k$ are associative and unitary 
$k$-algebras.
Given two such algebras, $A$ and $B$, the
category of 1-morphisms from $A$ to $B$ 
$$1\Hom_{\on{\mB im_k}}(A,B) = \on{\mathcal B\!\on{im}_k}(A,B)$$ is the
category of $A$-$B$-bimodules (where the left and right actions of $k$
are required to agree). Horizontal composition of 1-morphisms into and out
of $A$ is given by the tensor product over $A$.
Again, the degree zero part of the K\"unneth theorem for Hochschild
cohomology (c.f.\ \cite[X.7.4]{MacLane}) gives the isomorphism
\begin{eqnarray*}
  \mu\negmedspace:\ttr(M)\tensor\ttr(N)&\cong&
  \on{Center_A}(M)\tensor  \on{Center_A}(N)\\
  &\cong&\on{Center_A}(M\tensor N)\\
  &\cong&\ttr(M\tensor N).  
\end{eqnarray*}

This is closely related to 
Section \ref{tensor-algebra-Sec}:
Let $A$ and $B$ be associative and unitary finite dimensional
$k$-algebras. Let
$\on{\mathcal B\!\on{im}_k^f}(A,B)$
be the full subcategory of $\on{\mathcal B\!\on{im}_k}(A,B)$ whose
objects are bimodules that are finitely presented as right $B$-modules.
Then there is an equivalence of categories
\begin{eqnarray}
\label{eq:Funead}
  {\Funead\(\on{A-mod^f},\on{B-mod^f}\)}&\simeq& 
  {\on{\mathcal
      B\!\on{im}_k^f}(A,B)} \\
  \notag
  F&\longmapsto& F(A)\\
  \notag 
  {(-)\tensor_AM} &\longleftarrow& M.
\end{eqnarray}
Here ${\mathcal \Funead}$ stands for the category of
  $k$-linear right-exact functors and natural transformations between them. 
Under this equivalence, composition of functors
corresponds to 
composition of 1-morphisms in $\mathcal B\!\on{im}_k$.
\subsection{Bimodules, derived version}
\label{derived-bimodule-Sec}
The triangulated $k$-linear 2-ca\-te\-go\-ry $\mathcal
D\mathcal B\!\on{im}_k$ 
of \cite[2.4(b), 3.5]{Ganter:Kapranov} has the same objects
as $\mathcal B\!\on{im}_k$. The category
$$\on{1Hom}_{\mathcal D\mathcal B\!\on{im}_k}(A,B)$$
 is the derived
category of complexes of $A$-$B$-bimodules bounded
above. Horizontal composition is as above, but using the derived
tensor product.  
Let $A$ be a finite dimensional, associative and unitary
$k$-algebra, viewed
as an object of $\mathcal D\mathcal B\!\on{im}_k$. 
Let $K^\bullet$ be a 1-endomorphism of $A$. 
Then the categorical trace of $K^\bullet$ is the hypercohomology of $K^\bullet$,
\begin{eqnarray*}
  \ttr^\bullet(K^\bullet)& =& \Hom^\bullet_{\mathcal D^+\mathcal
    B\!\on{im}(A-A)}\(A,K^\bullet\) \\
  &=& \mathbb H^\bullet(A;K^\bullet).
\end{eqnarray*}
In the special case were $K^\bullet$ is a single $A$-$A$-bimodule $M$
situated in degree zero, 
\begin{eqnarray*}
  \ttr^\bullet(M)& =& \on{Ext}^\bullet_{A\tensor A^{\on{op}}}(A,M)\\
  &=& \on{HH}^\bullet(A,M)
\end{eqnarray*}
is the Hochschild cohomology of $M$.
Let $B$ be another finite dimensional, associative and unitary
$k$-algebra, and let $L^\bullet$ be a complex 
of $B$-$B$-bimodules, bounded above. 
Then
$$
\longmap{\mu}{\ttr^\bullet\(K^\bullet\)\tensor\ttr^\bullet\(L^\bullet\)}{\ttr\(K^\bullet\tensor 
  L^\bullet\)}
$$
 is the K\"unneth map.
In the special case where both $K^\bullet$ and $L^\bullet$ are single
bimodules situated in degree zero,
the K\"unneth theorem for Hochschild cohomology (c.f.\
\cite[X.7.4]{MacLane}) implies that $\mu$ is an isomorphism.
\subsection{Constructible sheaves}\label{sec:constr:mu}
The 2-category $\on{\mathbb R\mA n_\CC}$ of \cite[2.4.(d)]{Ganter:Kapranov} 
has as objects real analytic manifolds. For any two such manifolds
$X$ and $Y$, the category of 1-morphisms between them is defined as
$$\on{1Hom}_{\on{\mathbb R\mA n_\CC}}(X,Y) = \DCon(X\times Y),$$
the bounded derived category of ($\mathbb R$-)constructible sheaves of
$\CC$-vector spaces on $X\times Y$.\footnote{See \cite[Section
8.4]{Kashiwara:Schapira}, for background on constructible sheaves.} 
Horizontal composition of $$\mathcal K\in \on{1Hom}(X,Y)
\quad\text{and}\quad  
\mathcal \mL\in \on{1Hom}(Y,Z)$$ is given by the {\em derived
convolution}
$$
  \mL\circ\mK=R\pi_{13*}\(\pi_{12}\inv\mathcal
  K\tensor^L\pi_{23}\inv\mathcal L\).
$$
Here the $\pi_{ij}$ are the projections from $X\times Y\times Z$ 
to the according products with two factors, the tensor product is
taken over the constant sheaf $\ul\CC$, and the $\pi_{ij}\inv$ are the
respective inverse image sheaf functors.

Every object (``kernel'') $\mathcal K\in \on{1Hom}(X,Y)$ defines a functor
\begin{eqnarray*}
  F_{\mathcal K} \negmedspace :\DCon(X) &\longrightarrow &
  \DCon(Y)\\
  \mathcal F&\longmapsto & R\pi_{Y*}\(\pi_X\inv\mathcal
  F\tensor^L\mathcal K\).
\end{eqnarray*}
Let $\map gXX$ be analytic, and let 
$$\map\gamma{X}{X\times X}$$ be the inclusion of the graph
$\Gamma_g$. Consider the kernel 
$$
  \mathcal K_g = \ul\CC_{\Gamma_g} :=\gamma_*\ul\CC,
$$
the constant sheaf supported on $\Gamma_g$,
viewed as a complex situated in degree $0$.
This is a lift of $g_*$ to $\on{\RR\mA n_\CC}$. In other words, we have 
$
  F_{\mathcal K_g} = g_*.
$
We have\footnote{This is Proposition 5.5 in
  \cite{Ganter:Kapranov}. Note the sign mistake there:
  in the proof of
  Proposition 5.5, $\ul\CC_{X^g}\left[\on{codim}(X^g)\right]$ should be
  $\bigoplus\ul\CC_{X^g_\alpha}\left[-\on{codim}(X^g_\alpha)\right]$. 
  }
$$
  \ttr^\bullet\(\mathcal K_g\)\cong \bigoplus_\alpha H^{\bullet-codim
    X^g_\alpha}(X^g_\alpha;\ul\CC),
$$
where the sum runs over the connected components of the fixed point set.
Let $Y$ be a second object of $\on{\RR\mA n_\CC}$, and consider an analytic
automorphism $h$ of $Y$. 
Then the connected component
$X^g_\alpha\times Y^h_\beta$ of
$(X\times Y)^{(g,h)}$
has codimension 
\begin{equation}
  \label{eq:codim}
  \on{codim}(X^g_\alpha)+\on{codim}(Y^h_\beta)
\end{equation}
in $X\times Y$. 
The map
$$
  \longmap{\mu}{H^\bullet\(X^g,\ul\CC\)\tensor
    H^\bullet\(Y^h,\ul\CC\)}{H^\bullet\(X^g\times Y^h,\ul\CC\)} 
$$
is the K\"unneth isomorphism. 
Here the notation is short-hand for the direct sum over the connected
components of $X^g\times Y^h$, with the grading in each summand 
shifted by \eqref{eq:codim}.


%
\subsection{Coherent sheaves} 
\label{sec:DCoh-mu}
Another key example is the 2-category
$\on{\mV\!ar}_\CC$.\footnote{This is
  \cite[2.4.Exa(c)]{Ganter:Kapranov}. For a detailed account, see 
\cite{Caldararu:Willerton} or \cite{Huybrechts}.}
Its objects are smooth projective
varieties over $\CC$, and for two such varieties $X$ and $Y$, the
category of 1-morphisms from $X$ to $Y$ is 
$$
  \on{1Hom}(X,Y) = \DCoh(X\times Y),
$$  
the derived category of coherent sheaves on $X\times Y$. The formalism
is identical to that of the previous section, but now the derived
tensor product $\tensor^L$ is taken over the structure sheaf, and the
inverse image functors $\pi_{ij}\inv$ are replaced by the pull-back
functors $\pi_{ij}^*$. For an automorphism $g$ of $X$, the kernel
$$
  \mathcal K_g=\mathscr O_{\Gamma_g}:=\gamma_*\mathscr O_X
$$
lifts the auto-equivalence $g_*$ of $\DCoh(X)$ in the sense explained
in the previous section.
Its trace is identified by the following theorem:
\begin{Thm}
\label{thm:ttr(g)}
  Let $N^g$ be the normal bundle of
  the inclusion $$u\negmedspace : X^g\hookrightarrow X.$$ Then we have
  an isomorphism 
  $$
    \ttr^\bullet(g)\cong\bigoplus_\alpha
    \on{HH^{\bullet-codim(X^g_\alpha)}} \(X^g_\alpha,\det (N^g)\),
  $$ 
  where the direct sum is over the connected components 
  of the fixed point set $X^g$.
\end{Thm}
Again, 
$$
  \longmap{\mu}{H^\bullet\(X^g,\det(N^g)\)\tensor
    H^\bullet\(Y^h,\det(N^h)\)}{H^\bullet\(X^g\times
    Y^h,\det(N^{(g,h)})\)} 
$$
is the K\"unneth isomorphism. Again, the notation is short-hand for
the direct sum over the connected 
components of $X^g\times Y^h$, with the grading in each summand 
shifted by \eqref{eq:codim}.

\medskip
\begin{Pf}{of Theorem \ref{thm:ttr(g)}}
  We have a commuting diagram of closed immersions
  \begin{equation}
    \label{eq:Delta-X^g-diagram}
    \xymatrix{
    {X^g}\ar[0,2]^u \ar[2,0]_u\ar[1,1]^\delta && X\ar[2,0]^\gamma\\
    &X^g\times X^g\ar[1,1]^i&\\
    X \ar[0,2]_\Delta&& X\times X.
    }
  \end{equation}
  Here $\Delta$ and $\delta$ denote the inclusions of the respective
  diagonals, and $\gamma$ is the inclusion of the graph $\Gamma_g$.
  All squares in \eqref{eq:Delta-X^g-diagram} are cartesian. We
  have
  \begin{eqnarray*}
    \ttr^\bullet(g) &= & \on{Ext^\bullet_{X\times
        X}}\(\Delta_*\mathscr O_X,\gamma_*\mathscr O_X \) \\
    &= & \on{Ext^\bullet_{X}}\(L\gamma^*\Delta_*\mathscr O_X,\mathscr
    O_X \). 
  \end{eqnarray*}
  Write $\Theta$ for the composite of natural transformations
  $$
    \Theta\negmedspace : L\gamma^*\Delta_* \Longrightarrow
    L\gamma^*\Delta_* u_*Lu^* =
    L\gamma^*i_*\delta_*Lu^*\Longrightarrow u_*L\delta^*\delta_*Lu^*, 
  $$
  where the first map is given by the unit of the adjunction
  $Lu^*\dashv u_*$, and the second map is given by the base-change map
  of the upper right (triangle shaped) square in
  \eqref{eq:Delta-X^g-diagram}. 

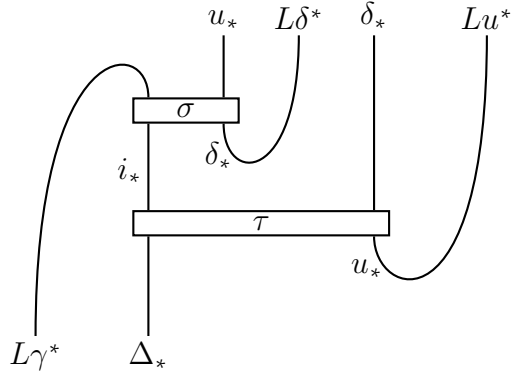
  \begin{tikzpicture}[inner sep=2pt,thick]
      \tikzset{
      short/.style={rectangle, draw=black, minimum width=1.4cm},
      long/.style={rectangle, draw=black, minimum width=3.4cm}
      }
      \node at (0,5) {};
      \node at (0,0) [name=Lgamma, anchor=north] {$L\gamma^*$};
      \node at (1.5,0) [name=Delta, anchor=north] {$\Delta_*$};
      \node at (3,1.5) [name=b, long] {$\tau$}; 
      \node at (2, 3) [name=a, short] {$\sigma$};
      \node at (2.5,4) [name=u, anchor=south] {$u_*$};
      \node at (3.5,4) [name=Ldelta, anchor=south] {$L\delta^*$};
      \node at (4.5,4) [name=delta, anchor=south] {$\delta_*$};
      \node at (6,4) [name=Lu, anchor=south] {$Lu^*$};
      \node at (7,2) [color=gray, anchor=west, text width=4.5cm] {
        The natural transformation $\Theta$ in string diagram
        notation. The squares labeled $\sigma$ and $\tau$ both
        denote the identity natural
        transformation of ${(\gamma u)_*}=(i\delta)_* = (\Delta
        u)_* $.};
      \node at (1.5,2.2) [anchor=east] {$i_*$};
      \node at (4.7,.9) [anchor=east] {$u_*$};
      \node at (2.7,2.4) [anchor=east] {$\delta_*$};
      
      \draw ([xshift=-1.5cm]b.south) -- (Delta.north);
      \draw ([xshift=-1.5cm]b.north) -- ([xshift=-.5cm]a.south);
      \draw (u.south) -- ([xshift=.5cm]a.north);
      \draw (Lgamma.north) .. controls (0,4) and (1.5,4)
      .. ([xshift=-.5cm]a.north);
      \draw (Ldelta.south) .. controls (3.5,2) and (2.5,2)
      .. ([xshift=.5cm]a.south);
      \draw([xshift=1.5cm]b.north) -- (delta.south);
      \draw([xshift=1.5cm]b.south) .. controls (4.5,.5) and (6,0) .. (Lu.south);
    \end{tikzpicture}

  We claim that $\Theta$ is an isomorphism. The claim implies the
  theorem, by the following calculation:
  \begin{eqnarray*}
    \on{Ext_X^\bullet}(L\gamma^*\Delta_*\mathscr O_X,\mathscr O_X) &\cong&
    \on{Ext_X^\bullet}(u_*L\delta^*\delta_*\mathscr O_{X^g},\mathscr
    O_X) \\
    &\cong&
    \on{Ext_{X^g\times X^g}^\bullet}(\delta_*\mathscr
    O_{X^g},\delta_*u^!\mathscr O_X) \\
    &=&
    \on{HH^\bullet}\(X^g,u^!\mathscr O_X\).
  \end{eqnarray*}
  Here $u^!$ is the right-adjoint of $u$. Explicitly, $u^!$ is given by 
  $$
    u^!\mathscr O_X\cong\bigoplus_{a}\det(N^g_\alpha)[-\on{codim}(X^g_\alpha)]
  $$
  (see e.g. \cite{Caldararu:derived}.)

  The proof that $\Theta$ is an isomorphism is an adaptation of the
  discussion in \cite[Appendix A]{Caldararu:Katz:Sharpe}. As the
  question is local, it will suffice to prove it over an open cover of
  $X\times X$. For $(x,y)\in X\times X\setminus \Delta(X^g)$, there
  exists an open neighourhood $U$ of $(x,y)$ with
  $U\cap\Delta=\emptyset$ or $U\cap\Gamma_g=\emptyset$. Since
  $\emptyset=\on{spec}(\{0\})$, the claim is trivially true over such
  a $U$.
  Let now $x\in X^g$. Writing 
  $$
    T_{\Delta(X)}:= \on{im}(T_\Delta),\quad
    T_{\Gamma_g}=\on{im}(T_\gamma),    
    \quad\text{and}\quad    T_{\Delta(X^g)}:=
    \on{im}(T_{i\delta}), 
  $$
  we have
  $$
    T_{\Delta(X)}\cap T_{X^g\times X^g}=T_{\Delta(X^g)} =
    T_{\Delta(X)}\cap T_{\Gamma_g}
  $$
  (inside $T_{X\times X}$).
  By Lemma \ref{lem:clean}, we have
  $$
    T_{\Delta(X)\cup (X^g\times X^g)} = T_{\Delta(X)}+ T_{X^g\times X^g}.
  $$
  We now apply Lemma \ref{lem:Y1} with $Z=X\times X$ and
  $W=\Delta(X)\cup (X^g\times X^g)$. Over a neighbourhood $U$ of $x$
  inside $X\times X$ this yields a refinement of
  \eqref{eq:Delta-X^g-diagram} looking as follows:
  $$
    \xymatrix{
      {X^g}\ar[0,2]^u \ar[2,0]_u\ar[1,1]^\delta && X\ar[2,0]^\gamma\\
      &X^g\times X^g\ar@{..>}[1,1]^i\ar[d]_j&\\
      X \ar[0,1]_{\Delta'}&Y\ar[r]_k& X\times X
    }
  $$
  (all objects intersected with $U$).
  Here $Y$ is a smooth subvariety of $X\times X$, the maps $j$, $k$,
  and $\Delta'$ are closed immersions, and 
  $$T_{Y,x} = T_{\Delta(X),x}+T_{X^g\times X^g,x}.$$
  The left solid square (with arrows $\delta$, $j$ $u$ and $\Delta'$) is
  still cartesian. Further, $T_{i\delta}$ defines an isomorphism
  $$
    T_{X^g,x} \cong T_{T,x}\cap T_{\Gamma_g,x}
  $$
  (intersection inside $T_{X\times X, x}$). Applying Lemma
  \ref{lem:Y2}, we conclude that (possibly after passing to a smaller
  neighbourhood of $x$) the right solid square
  (with arrows $u$, $\gamma,$, $j\delta$ and $k$) is also
  cartesian. Moreover, both intersections 
  $$
    X^g = \Delta(X)\cap_Y(X^g\times X^g)\quad\text{and}\quad
    X^g = Y\cap_{X\times X}\Gamma_g
  $$
  are clean and ``of the expected dimension'', meaning that 
  $$\dim(X^g)+\dim(Y)=\dim(X)+\dim(X^g\times X^g)$$
  and similarly for $Y\cap \Gamma_g$.
  Using \cite[p.V 20, Cor. to Theorem 4]{Serre}, we conclude that both
  these squares are {\em tor-indepent squares} in the sense of
  \cite[Def. (3.10.2)]{Lipman}.
  By \cite[Thm. (3.10.3)]{Lipman}, the base change
  tranformations 
  $$
    L\gamma^*k_*\Longrightarrow u_*L(j\delta)^*
  $$
  and
  $$
    L_j^*\Delta'_* \Longrightarrow \delta_*Lu^*
  $$
  are isomorphisms. Hence, over $U$ we have a natural isomorphism
  $$
    \Theta'\negmedspace : L\gamma^*\Delta_*\Longrightarrow
    u_*L\delta^*Lj^*\Delta'_* \Longrightarrow u_8L\delta^*\delta_*
    Lu^*. 
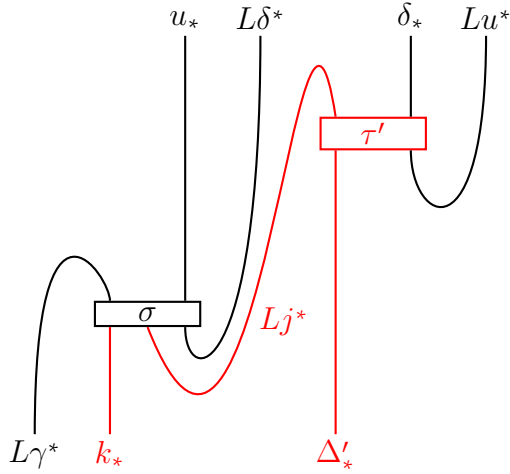
  $$

  \begin{tikzpicture}[inner sep=2pt,thick]
      \tikzset{
      short/.style={rectangle, draw=black, minimum width=1.4cm},
      long/.style={rectangle, draw=red, minimum width=1.4cm}
      }
      \node at (0,5) {};
      \node at (.5,-.3) [name=Lgamma, anchor=north] {$L\gamma^*$};
      \node at (1.5,-.3) [color=red, name=k, anchor=north] {$k_*$};
      \node at (4.5,-.3) [color=red, name=Delta', anchor=north] {$\Delta'_*$};
      \node at (5,3.7) [color=red, name=b, long] {$\tau'$}; 
      \node at (2, 1.3) [name=a, short] {$\sigma$};
      \node at (2.5,5) [name=u, anchor=south] {$u_*$};
      \node at (3.5,5) [name=Ldelta, anchor=south] {$L\delta^*$};
      \node at (5.5,5) [name=delta, anchor=south] {$\delta_*$};
      \node at (6.5,5) [name=Lu, anchor=south] {$Lu^*$};
      \node at (7.6,2.35) [color=gray, anchor=west, text width=3.8cm] {
        The natural transformation $\Theta'$ in string diagram
        notation. The squares labeled $\sigma$ and $\tau'$ 
        denote the identity natural
        transformations of ${(\gamma u)_*}=(kj\delta)_*$ and of
        $(\Delta'u)_*=(j\delta)_*$. The part 
        differing from the diagram for $\Theta$ is red.}; 
      \node at (3.4,1.2) [color=red, anchor=west] {$Lj^*$};
      
      \draw [color=red] ([xshift=-.5cm]b.south) -- (Delta'.north);
      \draw [color=red] ([xshift=-.5cm]a.south) -- (k.north);
      \draw [color=red] ([xshift=-.5cm]b.north) .. controls (4,7.2) and (3.5,-2.7) .. (a.south);
      \draw (u.south) -- ([xshift=.5cm]a.north);
      \draw (Lgamma.north) .. controls (.5,3.3) and (1.5,1.8)
      .. ([xshift=-.5cm]a.north);
      \draw (Ldelta.south) .. controls (3.5,.3) and (2.5,.3)
      .. ([xshift=.5cm]a.south);
      \draw([xshift=.5cm]b.north) -- (delta.south);
      \draw([xshift=.5cm]b.south) .. controls (5.5,2.5) and (6.5,2) .. (Lu.south);
    \end{tikzpicture}

  To see that $\Theta'=\Theta\vert_U$, we move $\tau'$ downward until
  it is 
  below $\sigma$, then
  straighten the red zig-zag line connecting $\sigma$ and $\tau'$ into a
  straight line labeled $j_*$, and finally
  note that $$k_*\tau'=\id_{(i\delta)_*}=\tau.$$
\end{Pf}
\section{Tensor products of 2-representations}
\label{sec:2reps}
Let $G$ be a finite group. Recall, e.g.\ from
\cite[Def.4.1.]{Ganter:Kapranov}, that a 2-representation of $G$ on a
$k$-linear category $\mV$ consists of the following data:
\begin{itemize}
\item for each $g\in G$ an endofunctor $\varrho(g)$ of $\mV$,
\item for any $g,h\in G$, a natural isomorphism
$$
  \phi_{g,h}\negmedspace :\varrho(g)\varrho(h)
  \stackrel\cong\Longrightarrow\varrho(gh), 
$$
and
\item a natural isomorphism $$
  \phi_1\negmedspace :\varrho(1)\stackrel\cong\Longrightarrow\id_\mV, 
$$
\end{itemize}
The $\phi_?$ are required to make associativity pentagons and unit
diagrams commute. In particular, for $g$, $h$ and $k$ in $G$ we have
an unambiguous double composition isomorphism
$$
  \phi_{g,h,k}\negmedspace :\varrho(g)\varrho(h)\varrho(k)
  \stackrel\cong\Longrightarrow\varrho(ghk).
$$ 
Following Bartlett\footnote{
See \cite[p.13]{Bartlett} for a translation of the axioms into string
diagram moves. Note that our string diagrams flow
upwards, while Bartlett's flow downwards.}, we use the string diagram
notation and write
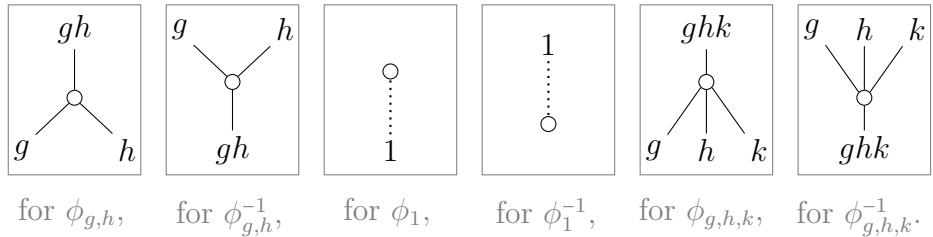
\begin{center}
\begin{tikzpicture}[inner sep=2pt, scale=.7]
   \tikzset{
   round/.style={circle, draw=black, radius=.01cm},
   }
  \node at (0,3) {};
  \node at (1,1) [name=phi, round] {};
  \node at (0,0) [name=g] {$g$};
  \node at (1,2.3) [name=gh] {$gh$};
  \node at (2,0) [name=h] {$h$};
  \node at (0.3,0.2) [name=g1] {};
  \node at (1,2.1) [name=gh1] {};
  \node at (1.7,0.2) [name=h1] {};
  \draw (gh.south) -- (phi.north);
  \draw (g.north east) -- (phi.south west);
  \draw (h.north west) -- (phi.south east);
  \node[draw=gray, fit=(g1) (gh1) (h1)] {};
  \node at (1,-.8) [gray, anchor=north, text width=1.5cm] {for $\phi_{g,h}$,};

  \node at (4,1.3) [name=phi, round] {};
  \node at (3,2.3) [name=g] {$g$};
  \node at (4,0) [name=gh] {$gh$};
  \node at (5,2.3) [name=h] {$h$};
  \node at (3.3,0.2) [name=g4] {};
  \node at (4,2.1) [name=gh4] {};
  \node at (4.7,0.2) [name=h4] {};
  \draw (gh.north) -- (phi.south);
  \draw (g.south east) -- (phi.north west);
  \draw (h.south west) -- (phi.north east);
  \node[draw=gray, fit=(g4) (gh4) (h4)] {};
  \node at (4,-.8) [gray, anchor=north, text width=1.5cm] {for $\phi_{g,h}\inv$,};

  \node at (7,1.5) [name=phi, round] {};
  \node at (7,0) [name=1] {$1$};
  \node at (6.3,0.2) [name=g2] {};
  \node at (7,2.1) [name=gh2] {};
  \node at (7.7,0.2) [name=h2] {};
  \draw [thick, dotted] (1.north) -- (phi.south);
 
  \node[draw=gray, fit=(g2) (gh2) (h2)] {};
  \node at (7,-.8) [gray, anchor=north, text width=1.25cm] {for $\phi_{1}$,};

  \node at (10,0.5) [name=phi, round] {};
  \node at (10,2) [name=1] {$1$};
  \node at (9.3,0.2) [name=g3] {};
  \node at (10,2.1) [name=gh3] {};
  \node at (10.7,0.2) [name=h3] {};
  \draw [thick, dotted] (1.south) -- (phi.north);
 
  \node[draw=gray, fit=(g3) (gh3) (h3)] {};
  \node at (10,-.8) [gray, anchor=north, text width=1.4cm] {for
    $\phi_{1}\inv$,};

  \node at (13,1.3) [name=phi, round] {};
  \node at (12,0) [name=g] {$g$};
  \node at (13,2.3) [name=ghk] {$ghk$};
  \node at (13,0) [name=h] {$h$};
  \node at (14,0) [name=k] {$k$};
  \node at (12.3,0.2) [name=g1] {};
  \node at (13,2.1) [name=gh1] {};
  \node at (13.7,0.2) [name=h1] {};
  \draw (ghk.south) -- (phi.north);
  \draw (g.north east) -- (phi.south west);
  \draw (h.north) -- (phi.south);
  \draw (k.north west) -- (phi.south east);
  \node[draw=gray, fit=(g1) (gh1) (h1)] {};
  \node at (13,-.8) [gray, anchor=north, text width=1.8cm] {for $\phi_{g,h,k}$,};

  \node at (16,1) [name=phi, round] {};
  \node at (15,2.3) [name=g] {$g$};
  \node at (16,0) [name=ghk] {$ghk$};
  \node at (17,2.3) [name=k] {$k$};
  \node at (16,2.3) [name=h] {$h$};
  \node at (15.3,0.2) [name=g4] {};
  \node at (16,2.1) [name=gh4] {};
  \node at (16.7,0.2) [name=h4] {};
  \draw (ghk.north) -- (phi.south);
  \draw (g.south east) -- (phi.north west);
  \draw (k.south west) -- (phi.north east);
  \draw (h.south) -- (phi.north);
  \node[draw=gray, fit=(g4) (gh4) (h4)] {};
  \node at (16,-.8) [gray, anchor=north, text width=1.8cm] {for $\phi_{g,h,k}\inv$.};
\end{tikzpicture}
\end{center}
We recall further that $\varrho(g\inv)$ is a weak inverse, and hence a
two-sided adjoint, of $\varrho(g)$. 
For the adjunction $\varrho(g\inv)\dashv\varrho(g)$,  
we have
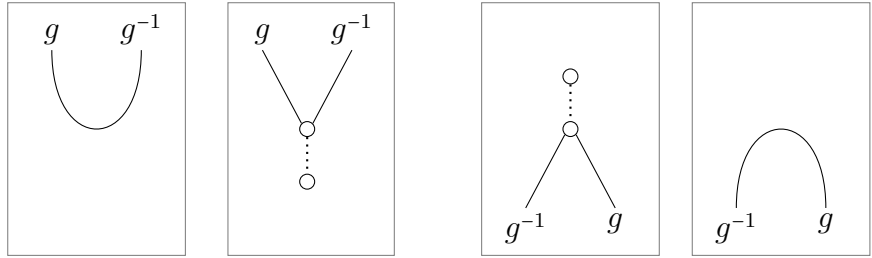
\begin{center}
  \begin{tikzpicture}[inner sep=2pt, scale=.7]
    \tikzset{
      round/.style={circle, draw=black, radius=.01cm},
      }
      \node at (0.65,0) [name=g, anchor=south] {$g$};
      \node at (2.35,0) [name=f, anchor=south] {$g\inv$};
       \node at (0.5,0) [name=a] {};
      \node at (2.5,-3) [name=b] {};

      \draw (g.south) .. controls (0.65,-2) and (2.35,-2) .. (f.south);
      \node[draw=gray, fit= (b) (a)] {};  

      \node at (4.65,0) [name=g, anchor=south] {$g$};
      \node at (6.35,0) [name=f, anchor=south] {$g\inv$};
      \node at (5.5,-1.5) [name=phi, round] {};
      \node at (5.5,-2.5) [name=phi1, round] {};
      \node at (4.65,0) [name=a] {};
      \node at (6.5,-3) [name=b] {};

      \draw (g.south) -- (phi.north west);
      \draw (f.south) -- (phi.north east);
      \draw [thick, dotted] (phi.south) -- (phi1.north);
      \node[draw=gray, fit= (b) (a)] {};

      \node at (3.5,-4.2) [gray, anchor=north, text width=5.3cm]
      {The unit is given by
       $(\phi_{g,g\inv}\phi_1)\inv$ and denoted by either of these
       two diagrams,}; 

      \node at (9.65,-3) [name=f, anchor=north] {$g\inv$};
      \node at (11.35,-3) [name=g, anchor=north] {$g$};
      \node at (10.5,-1.5) [name=phi, round] {};
      \node at (10.5,-.5) [name=phi1, round] {};
      \node at (9.5,0) [name=a] {};
      \node at (11.5,-3) [name=b] {};

      \draw (g.north) -- (phi.south east);
      \draw (f.north) -- (phi.south west);
      \draw [thick, dotted] (phi1.south) -- (phi.north);
      \node[draw=gray, fit= (b) (a)] {};

      \node at (13.65,-3) [name=f, anchor=north] {$g\inv$};
      \node at (15.35,-3) [name=g, anchor=north] {$g$};
      \node at (13.5,0) [name=a] {};
      \node at (15.5,-3) [name=b] {};

      \draw (f.north) .. controls (13.65,-1) and (15.35,-1) .. (g.north);
      \node[draw=gray, fit= (b) (a)] {};  
      
      \node at (12.5,-4.2) [gray, anchor=north, text width=5cm]
      {and the counit is given by 
       $\phi_{g\inv,g}\phi_1$ and denoted by either of these
       two diagrams.}; 
    \end{tikzpicture}
  \end{center}
The axioms for the $\phi_{g,h}$ imply that these adjunctions compose
as one would hope for, i.e., the adjunction for $g$ composed with that
for $s$ gives the adjunction for $sg$. We refer the reader to
\cite{Bartlett}, in particular Lemma 6 (iv), 
for a detailed account of these facts and their expression in terms of
string diagram moves.

The categorical character of $\varrho$ consists of the data
$$
  X_\varrho(g) = \ttr(\varrho(g)), \quad g\in G
$$
and for each pair $(s,g)$ the isomorphism
$$
  \longmap{\psi_s}{X_\varrho(g)}{X_\varrho(sgs\inv)},
$$
sending $\xi\in\ttr(\varrho(g))$ to 
\begin{center}
  \begin{tikzpicture}[inner sep=2pt, scale=.7]
    \tikzset{
      short/.style={rectangle, draw=black},
      round/.style={circle, draw=black, radius=.01cm},
      }
      \node at (0,0) [name=phi, round] {};
      \node at (0,1) [name=sgs] {$sgs\inv$};
      \node at (0,-2) [name=xi, short] {$\xi$};
      \node at (-1.4,-1.5) [name=s, anchor=east] {$s$};     
      \node at (0,-1.5) [name=g, anchor=south west] {$g$};     
      \node at (1.4,-1.5) [name=t, anchor=west] {$s\inv$};     
      \node at (0,-2.4) [name=a] {};
      \node at (0,.8) [name=b] {};
      \node at (0,1.9) {};

      \draw (phi.south west) .. controls (-5,-4) and (5,-4)
      .. (phi.south east);
      \node[draw=gray, fit= (b) (s) (t) (a)] {};  
      \draw (phi.south) -- (xi.north);
      \draw (sgs.south) -- (phi.north);
      \node at (0.2,-3.5) [gray, anchor=north, text width=3.9cm]
      {$\psi_s(\xi)\in\ttr(\varrho(sgs\inv))$}; 
    \end{tikzpicture}
  \end{center}
\label{page:psi}
and satisfying $\psi_1 = \id$ and $\psi_s\psi_t=\psi_{st}$
(see \cite[Prop.4.10]{Ganter:Kapranov} or \cite[Prop.16]{Bartlett} for
details). 
\begin{Def}
  Let $\varrho$ be a 2-representation of $G$ on a linear category $\mV$,
  and let $\pi$ be a 2-represenation of $H$ on a linear category
  $\mW$. We define the 2-representation $\varrho\boxtimes \pi$ on
  $\mV\boxtimes\mW$ by applying $-\boxtimes-$ to all the data of
  2-representation. 
\end{Def}
\begin{Cor}
  In the situation of Theorem \ref{tensor-Thm}, we have isomorphisms
  $$
    \isomap{\mu}{X_\varrho(g)\tensor X_\pi(h)}{X_{\varrho\boxtimes\pi}(g,h)}
  $$
  identifying $\psi^\varrho_s\tensor\psi^\pi_t$ with
  $\psi^{\varrho\boxtimes \pi}_{(s,t)}$.
  In orther words, $X_\varrho\tensor X_\pi$ and $X_{\varrho\boxtimes\pi}$
  are isomorphic as representations of the inertia groupoid
  $\Lambda(G\times H)$.
\end{Cor}
\begin{Pf}{}
  This follows from Theorem \ref{tensor-Thm} and inspection of the
  definition of $\psi_s$ in the proof of \cite[Prop.4.10.]{Ganter:Kapranov}.
\end{Pf}{}

A similar statement holds for $\widehat\boxtimes$ in the abelian case
and, more generally, in any monoidal 2-category where $\mu$ is an
isomorphism. 
\section{Categories of equivariant objects}
\label{V^G-Sec}
Let $\map\varrho GGL(V)$ be a complex representation of a finite group.
We write $V^G$ for the maximal $G$-invariant summand of $V$. 
The dimension of $V^G$ is
\begin{eqnarray}
  \label{inner-product-with-1-Eqn}
  \notag
  \dim(V^G) &=&  \langle V,\CC\rangle_G\\
 & =& \GG G\sum_{g\in G}\chi_\varrho(g),  
\end{eqnarray}
where $\chi_\varrho$ is the character of $\varrho$.

%
This section studies categorical analogues of $V^G$. These are the category
of equivariant objects and its generalizations.
\subsection{Background and Definitions}
The notion of {\em category of equivariant objects} already plays a role in
{\cite{Grothendieck:Tohoku}, see page 196 of [loc.cit].
\begin{Def}
Let $\mV$ be a linear category, acted upon from the left by a finite
group $G$. 
An {\em equivariant object} of $\mathcal{V}$ consists of an object
$x\in\ob(\mV)$ 
and a system of isomorphisms
$$\left\{ \map{ \epsilon_g}{x}{gx}\right\}_{ g\in G}$$
such that for any $g,h\in G$ the diagram
  \begin{equation}
    \label{eq:equivariant}
    \xymatrix{
    x\ar[r]^{\epsilon_g}\ar[d]_{\epsilon_{gh}} &
    {gx} \ar[d]^{g\epsilon_h} \\
    {(gh)x} &
    {g(hx)}\ar[l]^{\phi_{g,h,x}}
    }  
  \end{equation}
  is commutative.  
  As a consequence, we obtain the unit condition
  $$
    \epsilon_1 = \phi_{1,x}^{-1}: x\to 1x.
  $$
  A {\em map of
  equivariant objects} 
  $$
    \map f{(x,\epsilon)}{(y,\delta)}
  $$
  is a morphism $f\in\Hom(x,y)$ satisfing $\forall g\in G$
  $$
    g(f)\circ\epsilon_g = \delta_g \circ f. 
  $$
  The {category of equivariant objects} and their maps is denoted $\mV^G$.
\end{Def}
Let $X$ be a $G$-space. Then the categories of
equivariant $k$-vector bundles over $X$, equivariant
sheaves over $X$, and equivariant spaces over 
$X$ have interpretations as categories of equivariant objects. 
Below, we will discuss some of these examples in
detail. Maps of 2-representations form another example:
\begin{Def}\label{def:1Hom_G}
  Let $\mC$ be a 2-category, and let $\pi$ and $\varrho$ be
  2-re\-pre\-sen\-ta\-tions of $G$ on objects $W$ and $V$ of $\mC$.
  Then $G\times G$ 
  acts on the category $\on{1Hom}(W,V)$ by
  $$
    F\longmapsto \varrho(g)F\!\medspace\pi(s\inv).
  $$
  The coherence isomorphisms of this action are inherited from those
  of $\pi$ and $\varrho$. 

  View $G$ as a subgroup of $G\times G$ via the
  diagonal inclusion. We write
  $$
    \on{1Hom_G}(W,V):= \on{1Hom}(W,V)^G,
  $$
  for the {\em category of equivariant 1-homohorphisms and
    $G$-invariant 2-homomorphisms between them};
  and
  $$
    \on{1Hom_G^{triv}}(W,V):=
    \on{1Hom}(W,V)^{G\times G} 
  $$
  for the category of {\em trivialized $G$-equivariant 1-morphisms and
  tri\-via\-li\-za\-tion-preserving $G$-invariant 2-homomorphisms}.
Explicitly, an object of $\on{1Hom_G}(W,V)$ consists of a 1-morphism
$\map fWV$ together with a compatible family of 2-isomorphisms
$$
  \eta_g\negmedspace : F \Longrightarrow \varrho(g) F \pi(g).
$$
An alternative point of view is to consider the flip 2-isomorphisms
$\tau_g$ given by the string diagrams
\begin{center}
  \begin{tikzpicture}[inner sep=2pt, scale=.7]
    \tikzset{
      long/.style={rectangle, draw=black, minimum width=2cm},
      }
      \node at (0,0) [name=F1] {$F$};
      \node at (0,3.4) [name=F2, anchor=south] {$F$};
      \node at (0,1.7) [name=eta, long] {$\eta_{g\inv}$};
      \node at (1.2,3.4) [name=pig, anchor=south] {$\pi(g)$};
      \node at (-2.4,0) [name=rhog] {$\varrho(g)$};
      \node at (1.2,3.4) [name=a] {};
      \node at (-2.4,0.2) [name=b] {};
      \node at (0,4.5) [] {};

      \draw (rhog.north) .. controls (-2.4,3.8) and (-1.2,3.2)
      .. ([xshift=-1.2cm]eta.north);
      \draw (F2.south) -- (eta.north);
      \draw (pig.south) -- ([xshift=1.2cm]eta.north);
      \draw (F1.north) -- (eta.south);
      \node[draw=gray, fit= (a) (b)] {};  
      \node at (0.2,-1.2) [gray, anchor=north, text width=3.9cm]
      {Definition of $\tau_g$};

      \node at (5.6,0) [blue, name=F1] {$F$};
      \node at (5.6,3.4) [blue, name=F2, anchor=south] {$F$};
      \node at (5.6,1.7) [name=tau, blue] {$\bullet$}; 
      \node at (6.6,3.4) [name=pig, anchor=south] {$g$};
      \node at (4.6,0) [name=rhog] {$g$};
      \node at (6.6,3.4) [name=a] {};
      \node at (4.6,0.2) [name=b] {};

      \draw [blue, thick] (F1.north) -- (F2.south);
      \draw (rhog.north) .. controls (4.6,1.4) and (6.6,2) ..(pig.south);
      \node[draw=gray, fit= (a) (b)] {};  
      \node at (6.6,-1.2) [gray, anchor=north, text width=3.9cm]
      {Notation for $\tau_g$}; 

      \node at (10.6,0) [blue, name=F1] {$F$};
      \node at (10.6,3.4) [blue, name=F2, anchor=south] {$F$};
      \node at (10.6,1.7) [name=tau, blue] {$\bullet$}; 
      \node at (9.6,3.4) [name=pig, anchor=south] {$g$};
      \node at (11.6,0) [name=rhog] {$g$};
      \node at (9.6,3.4) [name=a] {};
      \node at (11.6,0.2) [name=b] {};

      \draw [blue, thick] (F1.north) -- (F2.south);
      \draw (rhog.north) .. controls (11.6,1.4) and (9.6,2) ..(pig.south);
      \node[draw=gray, fit= (a) (b)] {};  
      \node at (11.6,-1.2) [gray, anchor=north, text width=3.9cm]
      {Notation for $\tau_g\inv$}; 
    \end{tikzpicture}
  \end{center}
and satisfying
$$
  \tau_1 = \phi_1F\phi_1\inv
$$
(horizontal composition but vertical inverse) and
$$
  \tau_{gh}\circ\(\phi_{g,h}F\) =
  \(F\phi_{g,h}\)\circ\(\tau_gh\)\circ\(g\tau_h\) 
$$
(see \cite[(5)]{Bartlett}).
The data $\(F,\{\tau_g\}\)$ are eqivalent to those of
$\(F,\{\eta_g\}\)$.
An object of $\on{1Hom^{triv}_G}$ can then be thought of as a
$G$-equivariant 1-morphism $\(F,\{\tau_g\}\)$ from $V$ to $W$ together
with a family of {\em (left) trivilization 2-isomorphisms}
$$
  \epsilon_g\negmedspace : F\Longrightarrow \varrho(g)F,
$$
satisfying \eqref{eq:equivariant} 
and compatible with $\{\tau_s\}$ in the following sense:
\begin{center}
  \begin{tikzpicture}
      \node at (5.0,0.2) [blue, name=F1] {$F$};
      \node at (5.0,2) [blue, name=F2, anchor=south] {$F$};
      \node at (5.0,1) [name=eps, red] {$\bullet$}; 
      \node at (4.4,2) [name=pig, anchor=south] {$g$};
      \node at (6.0,2) [name=a] {};
      \node at (5.0,0.2) [name=b] {};

      \draw [blue, thick] (F1.north) -- ([yshift=-.06cm]eps.center);
      \draw [blue, thick] ([yshift=.08cm]eps.center) -- (F2.south);
      \draw ([xshift=-.08cm]eps.center) .. controls (4.3,1.3) and (4.45,1.9) ..(pig.south);
      \node[draw=gray, fit= (F1) (pig)] {};  
      \node at (5.2,-0.4) [gray, anchor=north, text width=3.9cm]
      {Notation for $\epsilon_g$};     
  \end{tikzpicture}
%
  \begin{tikzpicture}[inner sep=2pt];
      \tikzset{
      round/.style={circle, draw=black, radius=.01cm}}
      \node at (5.0,-.5) [blue, name=F1, anchor=north] {$F$};
      \node at (5.0,2.3) [blue, name=F2, anchor=south] {$F$};
      \node at (5.0,1) [name=eps, red] {$\bullet$}; 
      \node at (5.0,0.3) [name=tau, blue] {$\bullet$};
      \node at (4.2,1.7) [name=phi, round] {$$};  
      \node at (5.8,-.5) [name=rhos, anchor=north] {$s$};
      \node at (4.2,2.3) [name=gs, anchor=south] {$gs$};  
      
      \node at (4.35,1.15) [anchor=south west] {$g$};
       \node at (5,0.2) [name=b] {};

      \draw (phi.north) -- (gs.south);
      \draw [blue, thick] (F1.north) -- ([yshift=-.05cm]eps.center);
      \draw [blue, thick] ([yshift=.08cm]eps.center) -- (F2.south);
      \draw ([xshift=-.08cm]eps.center) .. controls (4.3,1.1) and
      (4.25,1.6) ..(phi.south east);
      \draw ([xshift=-.08cm]tau.center) .. controls (4.1,.6) and
      (4.25,1.6) ..(phi.south);
      \draw ([xshift=.08cm]tau.center) .. controls (5.6,0.2) and
      (5.7,.1) ..(rhos.north);
      \node[draw=gray, fit= (rhos) (gs)] {};

      \node at (8.0,-.5) [blue, name=F1,anchor=north] {$F$};
      \node at (8.0,2.3) [blue, name=F2, anchor=south] {$F$};
      \node at (8.0,1) [name=tau, blue] {$\bullet$}; 
      \node at (8.0,0.3) [name=eps, red] {$\bullet$};
      \node at (7.2,1.7) [name=phi, round] {};  
      \node at (8.8,-.5) [name=rhos, anchor=north] {$s$};
      \node at (7.2,2.3) [name=gs, anchor=south] {$gs$};  
      
      \node at (7.6,.7) [name= b, anchor=north east] {$s\inv gs$};

      \draw (phi.north) -- (gs.south);
      \draw [blue, thick] (F1.north) -- ([yshift=-.05cm]eps.center);
      \draw [blue, thick] ([yshift=.08cm]eps.center) -- (F2.south);
      \draw ([xshift=-.08cm]eps.center) .. controls (7.3,0.3) and
      (7.25,1.6) ..(phi.south);
      \draw ([xshift=-.08cm]tau.center) .. controls (7.1,1.3) and
      (7.36,1.65) ..(phi.south);
      \draw ([xshift=.08cm]tau.center) .. controls (8.8,.8) and
      (8.7,.1) ..(rhos.north);
      \node[draw=gray, fit= (rhos) (b) (gs)] {};  
      \node at (6.7,-1) [gray, anchor=north, text width=6cm]
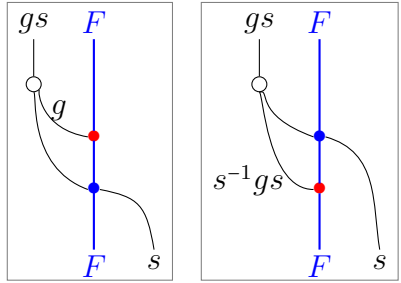
      {These are equal to each other.};     
  \end{tikzpicture}
\end{center}
Equivalently, one may work with the right trivialization isomorphisms
$$
  \epsilon^R_g := \tau_g\circ\epsilon_g
$$
with symmetric notation and compatibility requirements to those of
$\epsilon_g$ (here $\eta_{(g,1)}$ is replaced by $\eta_{(1,g\inv)}$).
In the special case where $F=\id_{W}$, we say that
the $G$-action on $W$ is trivialized (by $\{\epsilon_g\}$).
\end{Def}
\begin{Exa}\label{exa:composites}
  Let $\map FVW$ and $\map HWZ$ be equivariant
  1-morphisms, write $\{\tau_g\}$ and $\{\sigma_g\}$ for the
  respective families of twist maps. Then the composition $HF$,
  together with the twist maps $(H\tau_g)\circ(\sigma_gF)$, is again
  an equivariant 1-morphism. The horizontal composition of
  $G$-invariant 2-morphisms remains 
  $G$-invariant. 

  If $H$ is trivialized by $\{\epsilon_s\}$ then $HF$
  inherits the trivialization maps $\epsilon_gF$.
  If $F$ is trivialized with right trivialization maps
  $\{\delta_s^R\}$ then $HF$
  inherits a trivialization with right trivialization maps
  $\{H\delta_s^R\}$.  
  If both $H$ and $F$ are trivialized then the inherited trivilizations of
  $HF$ might not agree.
\end{Exa}
\begin{Exa}
  The identity 1-automorphism of $V$ is equivariant, with the $\eta_g$
  given by the units of the adjunctions
  $g\inv\dashv g$. The fact that these $\eta_g$
  compose in the expected way is \eqref{eq:equivariant}.

  If there exists a family $\{\epsilon_g\negmedspace:\id_V\Rightarrow
  \varrho(g)\}$ satisfying
  \eqref{eq:equivariant} then
  $\eta_g$ is identified with the horizontal composition 
  \begin{equation}
    \label{eq:epsgoepsginv}
    \eta_{g} = \epsilon_g\epsilon_{g\inv}.
  \end{equation}
  It follows that $\{\epsilon_g\}$ 
  is automatically compatible with $\{\tau_g\}$
  So, $\varrho$ is trivialized by $\{\epsilon_g\}$. 
  Further, $\epsilon_g^R=\epsilon_g$. 

  In this situation, we have 
  fully faithful embeddings
  \begin{eqnarray}\label{eq:faithful1}
    \on{1End}(V) & \longrightarrow & \on{1End_G}(V) \\
    \notag  F & \longmapsto & \(F,\{\epsilon_{g}F\epsilon_{g\inv}\}\).  
  \end{eqnarray}
  and
  \begin{eqnarray}\label{eq:faithful2}
    \on{1Hom_G}\(V,W\) & \longrightarrow &   \on{1Hom_G^{triv}}\(V,W\)\\
    \notag H & \longmapsto & \(H,\{H\epsilon_g\}\)
  \end{eqnarray}
  for any other object $W$ with $G$-action. 
\end{Exa}
We are now ready to formulate the universal property characterizing
$\mV^G$.
\begin{Prop}
\label{V^G-Prop}
  There is a $G$-action on $\mV^G$, which is canonically trivialized
  by the structure maps $\{\epsilon_g\}$. 
  The forgetful functor
  $$
    \longmap U{\mV^G}\mV
  $$
  is $G$-equivariant. Assume we are given another category $\mW$ with a
  $G$-action. Then composition with $U$ defines an equivalence of
  categories 
  \begin{eqnarray*}
    \on{\mF\!un}_G (\mW,\mV^G)
    &\stackrel \simeq\longrightarrow & \on{\mF\!un}_G^{triv} 
    (\mW,\mV).
  \end{eqnarray*}
  Here the trivialization of $UF$ is inherited from the trivialization
  $\{\epsilon_g\}$ of $\id_{\mathcal V^G}$ (see Example
  \ref{exa:composites}). In pictures: 
\end{Prop}
  \begin{tikzpicture}
      \node at (5.3,0.2) [blue, name=F1] {$F$};
      \node at (5.3,2) [blue, name=F2, anchor=south] {$F$};
      \node at (5.3,1) [name=eps, red] {$\bullet$}; 
      \node at (4.7,2) [name=pig, anchor=south] {$g$};
      \node at (6.3,2) [name=a] {};
      \node at (5.3,0.2) [name=b] {};

      \draw [blue, thick] (F1.north) -- ([yshift=-.06cm]eps.center);
      \draw [blue, thick] ([yshift=.08cm]eps.center) -- (F2.south);
      \draw ([xshift=-.08cm]eps.center) .. controls (4.6,1.3) and (4.75,1.9) ..(pig.south);
      \node[draw=gray, fit= (F1) (pig)] {};  

      \node at (7.25,0.2) [blue, name=F1] {$F'$};
      \node at (7.25,2) [blue, name=F2, anchor=south] {$F'$};
      \node at (6.75,0.2) [blue, name=U1] {$U$};
      \node at (6.75,2) [blue, name=U2, anchor=south] {$U$};
      \node at (7.0,1) [name=eps, red] {$\bullet$}; 
      \node at (6.4,2) [name=pig, anchor=south] {$g$};
      \node at (8.0,2) [name=a] {};
      \node at (7.0,0.2) [name=b] {};

      \draw [blue, thick] (F1.north) -- (F2.south);
      \draw [blue, thick] (U1.north) -- (U2.south);
      \draw ([xshift=-.08cm]eps.center) .. controls (6.3,1.3) and (6.45,1.9) ..(pig.south);
      \node at (6.75,1.1) [name=tau, blue] {$\bullet$}; 

      \node[draw=gray, fit= (F1) (pig)] {};  
      \node at (6.2,-0.4) [gray, anchor=north, text width=3.9cm]
      {These two diagrams agree.};     

      \node at (9.1,0.2) [blue, name=F1] {$F$};
      \node at (9.1,2) [blue, name=F2,anchor=south] {$F$};
      \node at (9.6,2) [name=pig,anchor=south] {$g$};
      \node at (8.6,0.2) [name=rhog] {$g$};
      \node at (10.1,2) [name=a] {};
      \node at (8.1,0.2) [name=b] {};

      \draw [blue, thick] (F1.north) -- (F2.south);
      \draw (rhog.north) -- (pig.south);

      \node at (9.1,1.23) [name=tau, blue] {$\bullet$}; 

      \node[draw=gray, fit= (pig) (rhog)] {};

      \node at (11,-0.4) [gray, anchor=north, text width=4.7cm]
      {These two diagrams agree.};       

      \node at (11.65,0.47) [blue, name=F1, anchor=north] {$F'$};
      \node at (11.65,2) [blue, name=F2,anchor=south] {$F'$};
      \node at (11.15,0.47) [blue, name=u1, anchor=north] {$U$};
      \node at (11.15,2) [blue, name=u2,anchor=south] {$U$};
      \node at (11.9,2) [name=pig,anchor=south west] {$g$};
      \node at (10.9,0.47) [name=rhog, anchor=north east] {$g$};
      \node at (12.4,2) [name=a] {};
      \node at (10.4,0.2) [name=b] {};

      \draw [blue, thick] (F1.north) -- (F2.south);
      \draw [blue, thick] (u1.north) -- (u2.south);
      \draw (rhog.north) -- (pig.south);

      \node at (11.65,1.49) [name=tau, blue] {$\bullet$}; 
      \node at (11.15,.97) [name=tau, blue] {$\bullet$}; 

      \node[draw=gray, fit= (pig) (rhog)] {};  

      \node at (13.4,0.2) [blue, name=F1] {$F$};
      \node at (13.4,2) [blue, name=F2, anchor=south] {$F$};
      \node at (13.4,1) [name=eps, red] {$\bullet$}; 
      \node at (14,2) [name=pig, anchor=south] {$g$};
      \node at (14.4,2) [name=a] {};
      \node at (13.4,0.2) [name=b] {};

      \draw [blue, thick] (F1.north) -- ([yshift=-.06cm]eps.center);
      \draw [blue, thick] ([yshift=.08cm]eps.center) -- (F2.south);
      \draw ([xshift=.08cm]eps.center) .. controls (14.1,1.3) and (13.95,1.9) ..(pig.south);
      \node[draw=gray, fit= (F1) (pig)] {};  
      \node at (15,-0.4) [gray, anchor=north, text width=3.5cm]
      {These two diagrams agree.};     

      \node at (15.5,0.2) [blue, name=F1] {$F'$};
      \node at (15.5,2) [blue, name=F2, anchor=south] {$F'$};
      \node at (15,0.2) [blue, name=U1] {$U$};
      \node at (15,2) [blue, name=U2, anchor=south] {$U$};
      \node at (15.25,1) [name=eps, red] {$\bullet$}; 
      \node at (15.85,2) [name=pig, anchor=south] {$g$};
      \node at (15.25,2) [name=a] {};
      \node at (14.85,0.2) [name=b] {};

      \draw [blue, thick] (F1.north) -- (F2.south);
      \draw [blue, thick] (U1.north) -- (U2.south);
      \draw ([xshift=.08cm]eps.center) .. controls (15.95,1.3) and (15.8,1.9) ..(pig.south);
      \node at (15.5,1.1) [name=tau, blue] {$\bullet$}; 

      \node[draw=gray, fit= (U1) (pig)] {};  
  \end{tikzpicture}

\begin{Pf}{}
  On objects of $\mV^G$, the $G$-action is defined by
  $$
    s\negmedspace:
    \{x,\epsilon_g\}\longmapsto\{sx,s\epsilon_{s\inv gs}\}. 
  $$ 
  All the other data of $G$-action are inherited from the action
  on $\mV$.
  Note that we are forced to define $\epsilon_g s$ in this manner in order for
  $\epsilon_g$ to be a natural transformation from $\id$ to $g$ in $\mV^G$:
  $$
    \(\epsilon_gs\)\circ\epsilon_s = \(g\epsilon_s\)\circ\epsilon_g
  $$
  is equivalent to
  \begin{equation}
    \label{eq:epsilon_gs}
    \epsilon_gs = \epsilon_{gs}\epsilon_s\inv = s\epsilon_{s\inv g s}.    
  \end{equation}
  By construction, the family $\{\epsilon_s\}$
  satisfies the coherence conditions \eqref{eq:equivariant}. 
  One checks that the functor of the proposition is a well-defined
  equivalence of categories. 
\end{Pf}{}

The universal property in Proposition \ref{V^G-Prop} should be compared that in
\cite[Prop. 4.4]{Ganter:Kapranov}. Although the latter looks less
general, it is equivalent to the one here: given
a trivialized functor $\map F\mW\mV$, we can apply 
\cite[Prop. 4.4]{Ganter:Kapranov} to the essential image of $F$ and 
obtain the result of Proposition \ref{V^G-Prop}. 
This argument does, however, not go through for representations in
general 2-categories, where we will work with the universal property
formulated here.
\begin{Def}
\label{def:V^G}
  Let $\mC$ be a 2-category, and let $\varrho$ be a 2-representation of
  $G$ on $V\in \ob(\mC)$.
  We will write 
  $$
    \longmap U{V^G}V
  $$
  whenever $V^G$ is an object of $\mC$ and $U$ is
  an equivariant 1-morphism from $V^G$ to $V$ such that $(V^G,U)$ 
  satisfies the universal property of Proposition \ref{V^G-Prop}.
\end{Def}
If it exists, the object $V^G$ is well defined up to trivialization
preserving equivariant 1-equivalence, which in turn is unique up to
canonical, trivialization preserving, $G$-invariant 2-morphism.
From now on, we fix $\longmap U{V^G}V$ as in Definition \ref{def:V^G}.
We write $\{\upsilon_g\}$ for the (left) trivialization maps of the
tivialization that $U$ inherits from $\id_{V^G}$, i.e., $\upsilon_g$
is the 2morphism $U\epsilon_g$ precomposed with the flip map $gU
\Rightarrow Ug$.

We will often
supress the representation $\varrho$ from the notation and write $g$
instead of
$\varrho(g)$. 
\begin{Prop}
\label{prop:adjoints}
  The 1-morphism $U$ possesses a left-adjoint $$A'\in
  1\Hom(V,V^G)$$ satisfying
  $$UA'=A:=\bigoplus_{g\in G}g.$$
\end{Prop}
\begin{Pf}{}
  The structure isomorphisms
  $$
    \isomap
    {\phi_{g,h,s\inv}}{\varrho(g)\varrho(h)\varrho(s\inv)}{\varrho(ghs\inv)}.
  $$
  make $A$ into a trivialized equivariant 1-morphism. 
  Write
  $
    \varphi_{s,A}\negmedspace : A\Rightarrow sA
  $
  for the trivialization 2-morphisms.
  Applying the universal property of $V^G$, we obtain a
  $G$-equivariant 1-morphism
  $$\longmap{A'}{V}{V^G}$$
  with $UA'=A$.
  Let $\iota$ be the composite
  $$
    \iota\negmedspace : \id_V\stackrel{\phi_1\inv}\Longrightarrow
    1\Longrightarrow A 
  $$
  where the second map is the inclusion of the $1$st summand.
  let
  $$
    \alpha\negmedspace : AU\Longrightarrow U
  $$
  be the 2-morphism
  $$
    \alpha:=\sum_{g\in G}\upsilon_g\inv.
  $$
  By \eqref{eq:equivariant}, we have
  $$
    \alpha\circ \(\phi_{s,A}U\) = \upsilon_s\circ\alpha.
  $$
  in other words, $\alpha$ is a trivialization preserving (equivariant) 
  2-morphism
  $$
    (AU,\{\varphi_{s,A}U\}) \Longrightarrow
    (U,\{\upsilon_s\}). 
  $$
  The universal property of $V^G$ yields a 2-morphism
  $$
    \alpha'\negmedspace : A'U \Longrightarrow 1
  $$
  with $U\alpha'=\alpha$.
  We need to show that $\iota$ and $\alpha'$ form the unit and counit
  of an adjunction $A'\rightfootline U$. 
  The identity
  $$
    \(U\alpha'\)\circ\(\iota U\)=\id_U 
  $$
  is immediate from the definitions and from $\epsilon_1
  = \phi_1\inv$. 
  By the universal property, it sufficies to check the second
  condition after horizontal composition with $U$. So, we are left to show
  $$
    \(\alpha A\)\circ\(A\iota\)=\id_{A}.    
  $$
  Restricted to the $g$th summand, this composite is the left-most in
  the following sequence of string diagrams: 
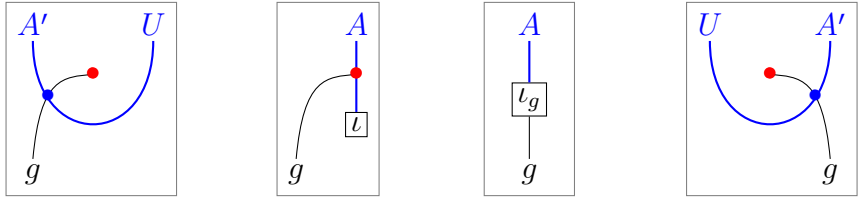
\begin{center}
\begin{tikzpicture}[inner sep=2pt]
  \node at (-3.8,0) [name=g0] {$g$};
  \node at (-3,1.333) [red, name=e0] {$\bullet$};
  \draw (g0.north) .. controls (-3.7,1.3) and (-3.4, 1.3) .. (e0.center);
  \node at (-3,1.333) [red, name=e0] {$\bullet$};
  \node at (-3.6,1.04) [blue, name=tau1] {$\bullet$};
  \node at (-2.2,2) [blue, name=u0] {$U$};
  \node at (-3.8,2) [blue, name=a'0] {$A'$};
  \draw [thick, blue] (u0.south) .. controls (-2.2,.3) and (-3.8,.3)
  .. (a'0.south);

  \node at (2.8,2) [blue, name=a1] {$A$};
  \node at (2.8,0) [name=g1] {$g$};
  \node at (2.8,1) [draw=black, rectangle, name=ig] {$\iota_g$}; 

  \draw (g1.north) -- (ig.south);
  \draw [thick, blue](a1.south) -- (ig.north);

  \node at (.5,2) [blue, name=a2] {$A$};
  \node at (.5,.666) [draw=black, rectangle, name=i] {$\iota$};
  \node at (-.3,0) [name=g2] {$g$};
  \draw [thick, blue] (a2.south) -- (i.north);
  \node at (.5,1.333) [red, name=e1] {$\bullet$};
  \draw (g2.north) .. controls (-.2,1.3) and (.1, 1.3) .. (e1.center);
  \node at (.5,1.333) [red, name=e1] {$\bullet$};

  \node at (6.8,0) [name=g3] {$g$};
  \node at (6,1.333) [red, name=e1] {$\bullet$};
  \draw (g3.north) .. controls (6.7,1.3) and (6.4, 1.3) .. (e1.center);
  \node at (6,1.333) [red, name=e1] {$\bullet$};
  \node at (6.6,1.04) [blue, name=tau1] {$\bullet$};
  \node at (5.2,2) [blue, name=u] {$U$};
  \node at (6.8,2) [blue, name=a'] {$A'$};
  \draw [thick, blue] (u.south) .. controls (5.2,.3) and (6.8,.3)
  .. (a'.south);
 
  \node[draw=gray, fit= (u0) (a'0) (g0)] {};    
  \node[draw=gray, fit= (a1) (g1), minimum width=1.2cm] {};    
  \node[draw=gray, fit= (a2) (g2)] {};    
  \node[draw=gray, fit= (u) (a') (g3)] {};    
  \node at (0,2.3) {};
\end{tikzpicture}  
\end{center}
Here $\iota_g$ denotes
  the inclusion of the $g$th summand.
The equality of the first three diagrams completes the proof. Their
equality to the fourth will be used below. 
\end{Pf}

\subsection{The twisted group algebra}
Let $G$ act on a linear category $\mathcal V$, and let $\mV^G$ be the
category of equivariant objects. We will think of this category as
resembling the category of
representations of $G$ ``in a twisted sense''. 
Representations of $G$ are the same as modules over the group algebra
$k[G]$, and it is natural to ask whether
there is an analog of $k[G]$ for this situation -- an associative algebra
acting on each equivariant object.
The answer is: yes, and very simple.

\medskip
Let $\varrho$ be a representation of $G$ on an object $V$ of a linear
2-category, and let $A$ be as in Proposition 
\ref{prop:adjoints}. 
Then, by general nonsense about adjunctions, $A=UA'$ is a monoid with
multiplication 
$$
  \alpha A'\negmedspace: A^2\Longrightarrow A
$$
and unit $\iota$.
Let $F$ be a 1-endomorphism of $V^G$. Then $M:= UFA'$ is a bimodule over
$A$ with the left- and right module structure given by ,
respecitively, $\alpha FA'$ and $UF\alpha'A'$.
\begin{Def}
  The {\em twisted group algebra} of $\varrho$
  is the space
  $$
    R_\varrho:=\ttr(A) = \bigoplus_{g\in G} \ttr (g),
  $$
  together with the $k$-algebra structure induced by the monoid
  structure of $A$.
\end{Def}
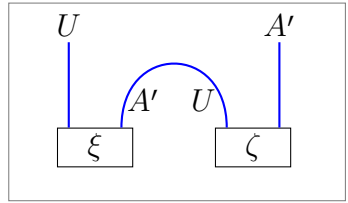
\begin{center}
\begin{tikzpicture}[inner sep=2pt, scale=.7]
   \tikzset{
   short/.style={rectangle, draw=black, minimum width=1cm},
   }
  \node at (0,1) [name=xi, short] {$\xi$};
  \node at (0.5,1.9) [anchor=west] {$A'$};
  \node at (2.4,1.9) [anchor=east] {$U$};
  \node at (-.5,3) [name=U, anchor=south] {$U$};
  \node at (-.5,3) [name=u] {};
  \node at (3,1) [name=zeta, short] {$\zeta$};
  \node at (3.5,3) [name=A, anchor=south] {$A'$};

  \draw [thick,blue] (U.south) -- ([xshift=-.5cm]xi.north);
  \draw [thick,blue] (A.south) -- ([xshift=.5cm]zeta.north);
  \draw [thick, blue] ([xshift=.5cm]xi.north) .. controls (.5,3) and (2.5,3)
  .. ([xshift=-.5cm]zeta.north); 
  \node[draw=gray, fit=(u) (zeta)] {};
  \node at (6,2) [gray, anchor=west, text width=4cm] {The product of
    elements $\xi$ and $\zeta$ of $\ttr(A)$};
\end{tikzpicture}
\end{center}
If $F$ is a 1-endomorphism of $V^G$, and $M=UFA'$ then $\ttr(M)$ is a
bimodule over $R_\varrho$.
\begin{Lem}
  Let $\xi$ be in the summand $\ttr(g)$ of $R_\varrho$, and let
  $\mu\in\ttr(M)$. Then their products are expressed by the following
  string diagrams:
  \begin{center}
  \begin{tikzpicture}[]
    \node at (1,2) [name=u] {$U$};
    \node at (1.5,1.3) [red,name=e] {$\bullet$};
    \node at (2,2) [name=F] {$F$};
    \node at (3,2) [name=a'] {$A'$};
    \node at (0,0) [draw=black, rectangle, name=xi] {$\xi$};
    \node at (2,0) [draw=black, rectangle, minimum width = 2.2cm,
    name=mu] {$\mu$};
    
    \draw [thick, blue] (u.south) -- ([xshift=-1cm]mu.north);
    \draw [thick, blue] (F.south) -- (mu.north);
    \draw [thick,blue] (a'.south) -- ([xshift=1cm]mu.north);
    \draw (xi.north) .. controls (.3,1) and (1,1.1) .. (e.center);
    \node at (1.5,1.3) [red,name=e] {$\bullet$};
    \node at (1,1.1) [blue] {$\bullet$};
    \node at (1.5,-1) [gray, text width=3.5cm] {Left multiplication};

    \node at (5,2) [name=u1] {$U$};
    \node at (6.5,1.3) [red,name=e1] {$\bullet$};
    \node at (6,2) [name=F1] {$F$};
    \node at (7,2) [name=a1] {$A'$};
    \node at (8,0) [draw=black, rectangle, name=xi1] {$\xi$};
    \node at (6,0) [draw=black, rectangle, minimum width = 2.2cm,
    name=mu1] {$\mu$};
    \node at (6.5,-1) [gray, text width=3.6cm] {Right multiplication};
   
    \draw [thick, blue] (u1.south) -- ([xshift=-1cm]mu1.north);
    \draw [thick, blue] (F1.south) -- (mu1.north);
    \draw [thick,blue] (a1.south) -- ([xshift=1cm]mu1.north);
    \draw (xi1.north) .. controls (7.7,1) and (7,1.1) .. (e1.center);
    \node at (6.5,1.3) [red,name=e1] {$\bullet$};
    \node at (7,1.1) [blue] {$\bullet$};

    \node [draw=gray, fit=(xi) (a')] {}; 
    \node [draw=gray, fit=(xi1) (u1)] {}; 
  \end{tikzpicture}    
  \end{center}
\end{Lem}
\begin{Pf}{}
  The diagram for the left multiplication is just the definition of
  $\alpha$. The diagram for the right multiplication follows from the
  equality of the third and the fourth diagram in the proof of
  Proposition \ref{prop:adjoints}.
\end{Pf}{}
If $F$ is of the form $\widetilde H$ for an equivariant 1-endomorphism
$(H,\{\tau_g\})$ of $V$ then $M=HA$, and
$$
  \ttr(HA)\cong \bigoplus_{g\in G}\ttr(Hg).
$$
\begin{Prop}
  The elements $\mu\in\ttr(Hs)$ and $\xi\in\ttr(g)$ multiply as
  follows:
  \begin{eqnarray*}
    \mu\cdot\xi &=& \(H\phi_{s,g}\)\circ\(\mu\xi\)\\
    \xi\cdot\mu &=& \(H\phi_{g,s}\)\circ\(\tau_gs\)\circ\(\xi\mu\).
  \end{eqnarray*}
  In particular, multiplication in $R_\varrho$ is horizontal
  composition followed by $\phi_{g,h}$
  \begin{eqnarray*}
    \ttr(\varrho(g))\tensor\ttr(\varrho(h))  &\longrightarrow & \ttr(\varrho(gh))
    \\
    \xi\tensor \zeta &\longmapsto & \phi_{g,h}\circ(\xi\zeta).
  \end{eqnarray*}
\end{Prop}
\begin{Lem}
\label{lem:j}
There is an isomorphism
of $k$-algebras
$$
  \isomap{j}{R_\varrho}{\on{2End}(U)},
$$
where the multiplication on
the target is vertical composition. For any 1-endomorphism $F$ of
$V^G$ there is an isomorphism of $R_\varrho$-bimodules
$$
  \isomap{j_F}{\on{2Hom}(\id_V,UFA')}{\on{2Hom}(U,UF).}
$$
\end{Lem}
\begin{Pf}{}
This is a consequence of the adjunction $A'\dashv U$. 
The isomorphism $j$ and its inverse are defined as
\begin{eqnarray*}
  j(\xi)&= &  \alpha\circ(\xi U)\\
  j\inv(\theta) &=& \(\theta A'\)\circ\iota,
\end{eqnarray*}
and similarly for $j_F$. More precisely,
$j_F$ is the 
second isomorphism on page 14 of \cite{Caldararu:Willerton} with
$\Psi=A'$, $\Phi = U$, 
  $\widehat\Theta=\id$, and $\widehat \Xi=UF$.
\end{Pf}{}

In the situation where $G$ acts on a linear category $\mV$, let
$(V,\{\epsilon_g\})$ be an equivariant object. Then we can
compose $j$ with the map 
\begin{eqnarray*}
  2\on{End}(U)&\longrightarrow & \on{End_\mV}(V)\\
  \theta&\longmapsto&\theta_V
\end{eqnarray*}
to obtain the promised action of $R_\varrho$ on $V$.
Explicitly,
$\xi = \(\xi_g\)_{g\in G}$
acts on $V$ by the endomorphism
$$
  \sum_g (\epsilon_g\inv \circ \xi_{g, V}).
$$

\begin{Exa}
Let $\mV=\on{Vect_k}$, and let the $G$-action be defined by means of a 2-cocycle
$c\in Z^2(G, k^*)$. So, we have the central extension
$$
  \map{p}{\tilde G}G
$$
with kernel $k^*$ and a bunch of 1-dimensional spaces
$$
  L_g = p^{-1}(g) \cup {0},
$$
as in \cite{Ganter:Kapranov}.
Then $\ttr (g) = L_g$,
and the algebra $R_\varrho$ is the direct sum of all the $L_g$. This is the standard
concept of the twisted group algebra associated to a central extension,
or the cocycle. An alternative definition of this algebra is by  a
basis $b_g$, for $g\in G$
with multiplication law
$$b_g\cdot b_h = c(g,h) \cdot b_{gh}.$$  
\end{Exa}
\subsection{The $G$-action on the twisted group algebra}
Let $F$ be an equivariant endofunctor of $V^G$. Then
the group $G$ acts on $\ttr(UFA')$ via the equivariance
isomorphisms of $\id_V$ and $UFA'$ (see Definition
\ref{def:1Hom_G}). In the situation where $F=\widetilde H$ for an
equivariant endofunctor $H$ of $V$, this action sends $\xi\in\ttr(Hg)$ to

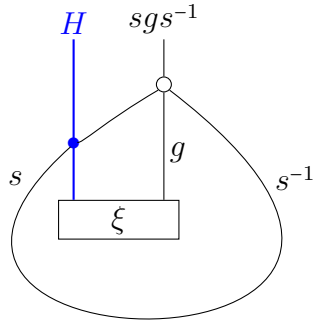
\begin{center}
\begin{tikzpicture}[inner sep=2pt]
      \tikzset{
      round/.style={circle, draw=black, radius=.01cm},
      long/.style={rectangle, draw=black, minimum width=1.6cm}
      }
      \node at (-.2,4.6) [name=g, anchor=south] {$sgs\inv$};      
      \node at (-1.4,4.6) [blue, name=H, anchor=south] {$H$};      
      \node at (0,5.5) [] {};      
      \node at (-0.2,4) [name=Phi, round] {};
      \node at (-.8,2.2) [name=xi, long] {$\xi$};
      \node at (-1.4,3.2) [name=tau, blue] {$\bullet$};
      \node at (-2,2.75) [anchor=east] {$s$};      
      \node at (1.2,2.75) [anchor=west] {$s\inv$};      
      \node at (-.2,3.1) [anchor=west] {$g$};      
      \node at (1.3,.5) [gray, anchor=north, text width=12cm] {The action of $G$ on
        $\ttr(HA)$ in string diagram notation};

      \draw (g.south) -- (Phi.north); 
      \draw (Phi.south west) .. controls (-1,3.5) and (-1.3,3.2) .. 
      ([xshift=.05cm, yshift=.03cm]tau.center);
      \draw ([xshift=-.04cm, yshift=-.03cm]tau.center) .. controls (-5,0) and (5,0) .. 
      (Phi.south east);
      \draw (Phi.south) -- ([xshift=.6cm]xi.north);
      \draw [thick, blue] (H.south) -- ([xshift=-.6cm]xi.north);
     
   \end{tikzpicture}   
\end{center}
\label{page:action}
In particular, the action on $\ttr(A)$ 
agrees with that of the isomorphisms
${\psi_s}$, defined on page \pageref{page:psi}. 

Consider the category $\on{Vect_G}(G)$ of $G$-equivariant vectorbundles on $G$
with respect to the conjugation action. 
The convolution product
$$
  (V\ostar W)_h := \bigoplus_{gs=h}V_g\tensor W_s 
$$
makes $\on{Vect_G}(G)$ a braided monoidal category, where the
braiding comes from the isomorphisms
$$
  V_g\tensor W_s\cong W_{s}\tensor V_{s\inv gs}
$$
induced by $\psi_s$.
A straight-forward manipulation of string diagrams\footnote{See
\cite[(2), (3), Lemma 6(i)]{Bartlett} for the relevant moves.} shows that 
$R_\varrho$ is a commutative algebra object of $\on{Vect_G}(G)$
with respect to this
(non-symmetric) braided monoidal structure.

\begin{Lem}
\label{lem:center}  
  The $G$-invariant part of the twisted group algebra $R_\varrho$ is
  contained in the center of $R_\varrho$. More generally, let $F$ be a
  1-endomorphism of $V^G$, viewed as an equivariant 1-endomorphsim via
  \eqref{eq:faithful1}. Then we
  have $$\ttr(M)^G\sub\on{Center}_{\ttr(A)}(\ttr(M)).$$  
\end{Lem}
\begin{Pf}{}
  Using the string diagram moves \cite[(2), (3), Lemma 6(i)]{Bartlett},
  one sees that 
  $\mu$ is in $\ttr(M)^G$ if and only if for each $s\in
  G$ the following two string diagrams are equal:
  
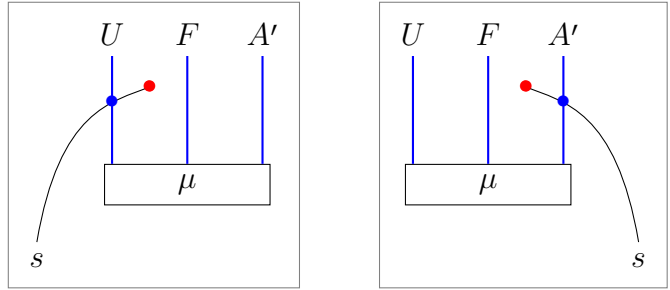
\begin{center}
  \begin{tikzpicture}[]
    \node at (1,2) [name=u] {$U$};
    \node at (1.5,1.3) [red,name=e] {$\bullet$};
    \node at (2,2) [name=F] {$F$};
    \node at (3,2) [name=a'] {$A'$};
    \node at (0,-1) [name=xi] {$s$};
    \node at (2,0) [draw=black, rectangle, minimum width = 2.2cm,
    name=mu] {$\mu$};
    
    \draw [thick, blue] (u.south) -- ([xshift=-1cm]mu.north);
    \draw [thick, blue] (F.south) -- (mu.north);
    \draw [thick,blue] (a'.south) -- ([xshift=1cm]mu.north);
    \draw (xi.north) .. controls (.3,1) and (1,1.1) .. (e.center);
    \node at (1.5,1.3) [red,name=e] {$\bullet$};
    \node at (1,1.1) [blue] {$\bullet$};

    \node at (5,2) [name=u1] {$U$};
    \node at (6.5,1.3) [red,name=e1] {$\bullet$};
    \node at (6,2) [name=F1] {$F$};
    \node at (7,2) [name=a1] {$A'$};
    \node at (8,-1) [name=xi1] {$s$};
    \node at (6,0) [draw=black, rectangle, minimum width = 2.2cm,
    name=mu1] {$\mu$};
   
    \draw [thick, blue] (u1.south) -- ([xshift=-1cm]mu1.north);
    \draw [thick, blue] (F1.south) -- (mu1.north);
    \draw [thick,blue] (a1.south) -- ([xshift=1cm]mu1.north);
    \draw (xi1.north) .. controls (7.7,1) and (7,1.1) .. (e1.center);
    \node at (6.5,1.3) [red,name=e1] {$\bullet$};
    \node at (7,1.1) [blue] {$\bullet$};

    \node [draw=gray, fit=(xi) (a')] {}; 
    \node [draw=gray, fit=(xi1) (u1)] {}; 
  \end{tikzpicture}    
  \end{center}
%
%
%
It follows that $\mu$ is in the center of $\ttr(M)$.
\end{Pf}
\begin{Lem}
\label{lem:jG}
  Let $F$ be an equivariant endofunctor of $V^G$. Then the isomorphism
  $j_{\widetilde H}$ of Lemma \ref{lem:j} induces an isomorphism of the
  $G$-invariant parts
  \begin{eqnarray*}
     \ttr(UFA')^G
     =2\Hom_G(\id_V,UFA')
    \cong 2\on{Hom_G}(U,UF).
  \end{eqnarray*}
\end{Lem}
\begin{Cor}
\label{cor:UHU}
  Let $H$ be an equivariant 1-endomorphism of $V$. Then we have an
  isomorphism
  $$
    2\on{End_G}(U,HU) \cong
    \(\bigoplus_{g\in G}\ttr(Hg)\)^G,
  $$
  where the $G$-action on the right-hand side is as on Page
  \pageref{page:action}. 
\end{Cor}

\subsection{Traces in $V^G$}
Throughout this section, we assume that we are given
a linear 2-representation of $G$ on $V$ such that $(V^G,U)$
exist. 
\begin{Def}
  Let $W$ be an object of a 2-category $\mC$. 
  Then the {\em dimension} or {\em center} of
  $W$ is 
  the categorical trace of the identity 1-morphism of $W$,
  $$
    \Cent(W) := \ttr(\id_W) = 2\Hom(\id_W,\id_W).
  $$
\end{Def}
Horizontal and vertical composition agree on $\Cent(W)$ and make it
into a commutative monoid.
If $\mC$ is 
$k$-linear then $\Cent(W)$ is a commutative $k$-algebra.
\begin{Thm}\label{thm:invariants}
  We have an isomorphism of $k$-algebras
  $$\isomap{i}{\Cent(V^G)}{\(\bigoplus_{g\in G}\ttr(g)\)^G}.$$
  For any 1-endomorphism $F$ of $V^G$ we have a bimodule isomorphism 
  $$
    \isomap{i_F}{\ttr(F)}{\ttr(UFA')^G}.
  $$
\end{Thm}

\begin{Pf}{}
  We have
  \begin{eqnarray*}
    \ttr(F) &= & 2\Hom\(\id_{V^G},F\)\\
    &= & 2\on{Hom_G}\(\id_{V^G},F\)\\
    &\cong & 2\on{Hom_G^{triv}}\(UF,UF\)\\
    &= & 2\on{Hom_G}\(UF,UF\)\\
    &\cong & \ttr(UFA')^G, 
  \end{eqnarray*}
  where the second equality is \eqref{eq:faithful1}, the 
  third isomorphism is the universal property of $(V^G,U)$, the
  next equality is \eqref{eq:faithful2}, and the last equality is 
  Lemma \ref{lem:jG}.
  If $F=\id_{V^G}$ we may apply Corollary \ref{cor:UHU} to identify the
  last term with
  $${\(\bigoplus_{g\in G}\ttr(g)\)^G}.$$

Explicitly, the isomorphism $i$ in the statement of the
theorem is given by
$$
  i\negmedspace :\zeta \longmapsto \(U\zeta A'\)\circ\iota.
$$
The fact that $i$ preserves the algebra multiplication is proved by
the equality of the following two string diagrams:

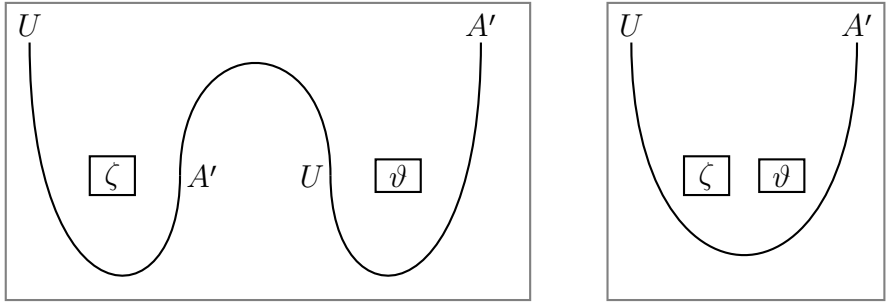
\begin{tikzpicture}[inner sep=2pt,thick]
      \tikzset{
      short/.style={rectangle, draw=black, minimum width=0.6cm},
      long/.style={rectangle, draw=black, minimum width=2.4cm}
      }
      \node at (0,4.5) {};
      \node at (-2.9,2) [name=zeta, short] {$\zeta$};
      \node at (0.9,2) [name=theta, short] {$\theta$};
      \node at (0,2) [anchor=east] {$U$};
      \node at (-2,2) [anchor=west] {$A'$};
      \node at (-2,2) [name=a] {};
      \node at (0,2) [name=u] {};
      \node at (-4,4) [name=U-up] {$U$};
     \node at (2.5,0.5) [name=c] {}; 
      \node at (2,4) [name=A] {$A'$};
      \node at (4,4) [name=U2] {$U$};
      \node at (7,4) [name=A2] {$A'$};
      \node at (5.5,0.5) [name=iota3] {};
      \node at (1.5,0) [gray] {The multiplication of $i(\zeta)$ and $i(\theta)$
        in $\ttr(A)$ equals $i(\zeta\theta)$.};
      \node at (5,2) [name=zeta, short] {$\zeta$};
      \node at (6,2) [name=theta, short] {$\theta$};
      
      \draw (U-up.south)  .. controls (-4,0) and (-2,0) .. (a.center);
      \draw (A.south)  .. controls (2,0) and (0,0) .. (u.center);
      \draw (a.center)  .. controls (-2,4) and (0,4) .. (u.center);
      \draw (U2.south) .. controls (4,0) and (7,0) .. (A2);
      \node[draw=gray, fit=(U-up) (c)] {};
     \node[draw=gray, fit=(U2) (iota3) (A2)] {};
    \end{tikzpicture}
\label{page:multiplicativity}

A similar identity of string diagrams shows that $i_F$ is a map of
bimodules. 
\end{Pf}
%

\begin{Cor}\label{2-inner-product-Cor}
Let $\varrho$ be a representation on an oject $V$ of a linear
2-category, and assume that $V^G$ exists. Then we have 
$$
  \dim\(\Cent(V^G)\) = \GG G\sum_{gh=hg}\chi_\varrho(g,h).
$$
\end{Cor}
\begin{Pf}{}
  We have
  \begin{eqnarray*}
    \dim\(\ttr(A)^G\) &=&\dim \bigoplus_{[g]\sub
  G}\ttr(g)^{C_g}\\
  &=& \sum_{[g]}\GG{C_g} \on{tr}(h\acts\ttr(g))\\
  &=& \GG G\sum_{gh=hg}\chi_\varrho(g,h).
  \end{eqnarray*}
\end{Pf}
\begin{Lem}
  Let $(H,\{\tau_g\})$ be an equivariant 1-endomorphism of $V$, and
  let $(\widetilde H,\{\widetilde \tau_g\})$ be such that $U\widetilde
  H$ and $HU$ agree as trivialized equivariant 1-morphisms.
  Then we have 
  $$
    \widetilde\tau_g \cong \epsilon_g \widetilde H\epsilon_g\inv
  $$ 
  (horizontal composition, vertical inverse). In other words,
  $\{\widetilde\tau_g\}$ agrees with the trivialization   
  \eqref{eq:faithful1}.
\end{Lem}
\begin{Cor}
  Let $(H,\{\eta_g\})$ is an equiariant 1-endomorphism of $V$ and let
  $(\widetilde H,\{\widetilde \eta_g\})$ then 
  $${\ttr(\widetilde H)}\cong{\(\bigoplus_{g\in G}\ttr(Hg)\)^G}.$$
\end{Cor}
\subsection{Inner products}
\begin{Def}
  Let $\varrho$ and $\pi$ be 2-representations of $G$ in a lax monoidal
  linear 2-category $(\mC,\boxtimes, {1})$,
  and assume that the object $(\varrho\boxtimes \pi)^G$ exists. Then we define the
  inner product of $\varrho$ and $\pi$ to be the $k$-vector space 
  $$\langle \varrho, \pi \rangle_G:= \Cent{} \((\varrho\boxtimes \pi)^G\).$$
\end{Def}
As an immediate consequence of the corollary, we obtain
$$
  \dim_k\langle\varrho,\pi\rangle_G = \GG G\sum_{gh=hg}\chi_\varrho(g,h)\cdot\chi_\pi(g,h).
$$
\section{Applications}
\label{applications-Sec}
\subsection{Projective representations}
Let $\mV=\on{Vect}_\CC$. Then a linear $G$-action on $\mV$ 
is classified by a $2$-cocycle
$$
  \map{c}{G\times G}{\CC},
$$
and $\mV^G$ is identified with the category of projective
representations of $G^{op}$ with central charge 
$$
  c^{op}(g,h) := c(h,g).
$$
These are pairs $$(W,\map{\varphi}{G^{op}}{\on{Aut}(W)}),$$ where $W$ is a
$k$-vector space, and $\varphi$ is a map satisfying
$$
  \varphi(gh) = c(g,h)\cdot \varphi(h)\circ\varphi(g).
$$
In particular, $\mV^G$ is a 2-vectorspace, and 
$$
  \dim\Cent(\mV^G)
$$
is the number of isomorphism classes of irreducible projective
representations of $G^{op}$ with central charge $c^{op}$.

On the other hand, 
$$
  \ttr(A)^G\cong\bigoplus_{[g]}\CC^{C_g},
$$
where the sum is over the conjugacy classes of $G$. 
The action of $h\in C_g$ on $\CC\cong \ttr(g)$ is
multiplication with
$\chi_{\varrho_c}(g,h)$.
\begin{Lem}
We have
$$
  \chi_{\varrho_c}(g,h) = \frac{c(h,g)}{c(g,h)}.
$$  
\end{Lem}
\begin{Pf}{}
Apply the cocycle condition to
$\delta(g,h,h\inv)$ to get
$$
  c(gh,h\inv)c(h,h\inv)\inv = c(g,h)\inv c(g,1).
$$
Substituting this into the formula of
\cite[Prop.5.1]{Ganter:Kapranov}), 
we obtain
\begin{eqnarray*}
    \chi_{\varrho_c}(g,h) &=&
     c(h,g) c(g,h)\inv c(g,1)c(1,1)\inv\\
  &=& c(h,g)c(g,h)\inv.
\end{eqnarray*}  
\end{Pf}{}
A conjugacy class $[g]\sub G$ is called {\em $c$-regular} if
  $$
    \(\forall h\in C_g\)\quad\(c(g,h)=c(h,g)\). 
  $$  
Note $[g]_G$ is $c$-regular if and only if $[g]_{G^{op}}$ is $c^{op}$-regular.
Hence Theorem \ref{thm:invariants} specializes to the following
result of Schur's:
\begin{Thm}[Schur]
  The number of isomorphisms classes of
  irreducible projective representations with multiplier $c$
  equals the number of $c$-regular conjugacy classes of
  $G$. 
\end{Thm}
\subsection{Algebras}
Let $A$ be an associative and unitary, finite dimensional $k$-algebra.
Assume that $G$ acts on $A$ from the left by (unit preserving) algebra
automorphisms. Let $\mV=\on{A-mod^f}$ be as in Section
\ref{tensor-algebra-Sec}. Then $G$ acts on $\mV$ from the right via
the right-exact functors
$$
  M\longmapsto Mg.
$$
Here $Mg$ is the right $A$-module 
that is isomorphic to $M$ as a $k$-vector space and whose $A$-action
is twisted by $g$ as follows:
$$
  (mg)\cdot a = (m\cdot g(a))g.
$$
Write $A\rtimes G$ for the crossed product algebra. So,
$$
  A\rtimes G = \bigoplus_{g\in G}Ag
$$
with multiplication
$$
  (a_1g_1)\cdot(a_2g_2) = (a_1\cdot g_1(a_2))(g_1g_2).
$$
The category of (right-)equivariant objects in $\mV$ is
$$
  \mV^G\simeq\on{(A\rtimes G)-mod^f}.
$$
The twisted group algebra $R_\varrho$ equals $$\bigoplus_{g\in
  G}\on{Center}_A(Ag).$$  
Hence Theorem \ref{thm:invariants} becomes the well known formula
\begin{equation}
  \label{eq:center-AG}
  \on{Center}(A\rtimes G)\cong  \bigoplus_{[g]}
  \(\on{Center}_A(Ag)\)^{C_g}.   
\end{equation}

\subsection{Bimodules}
Let $A$ be an associative and unitary $k$-algebra, acted upon by $G$ 
as in the previous section.
View $A$ as object of the 2-category
$\mathcal B\!\on{im}_k$ of Section \ref{sec:bimodules,strict}.
Let $U=A\rtimes G$, viewed as ($G\times G$-equivariant) $A\rtimes
G$-$A$-bimodule. Under the equivalence \eqref{eq:Funead} $U$
corresponds to the forgetful functor sending an equivariant $A$-module
to the underlying $A$-module. We would like to argue that the pair
$(A\rtimes G, U)$ satisfies the universal property of Definition
\ref{def:V^G}. This is not entirely true, but holds for a sufficiently
large class of test-objects: let $B$ be a second associative and
unitary $k$-algebra on which $G$ acts by algebra
automorphisms.\footnote{For the full universal property to by
  satisfied, we would need to allow $G$ to act by 1-automorphisms.}
Then composition with $U$ gives an equivalence of categories
$$\mathcal B\!\on{im}_{k,G\times G}(B,A)\simeq \mathcal B\!\on{im}_{k,
  G}(B,A\rtimes G).$$
This is sufficient to ensure that the proof of Theorem
\ref{thm:invariants} goes through, again yielding \eqref{eq:center-AG}.

We now view $A$ and $A\rtimes G$ as objects of the 2-category
$\mD\mB\!\on{im_k}$ (see Section \ref{derived-bimodule-Sec}), and we
write $\mathcal U$ for the complex consisting of $U$ situated in
degree zero and zeros elsewhere. 
For a test-object $B$ as above, we need to distinguish between the
derived category of equivariant bimodules and the category of
equivariant 1-morphisms: in general, the forgetful functor
$$
  \mathcal D\mathcal B\!\on{im}_{k,G}(B,A)\longrightarrow   \mathcal
  D\mathcal B\!\on{im}(B,A)^G = \on{1-Hom_G}(B,A), 
$$ 
is not an equivalence of categories. 

In addition to restricting our pool of test-objects as above, 
we need to modify the universal property in Definition
\ref{def:V^G}, replacing $\on{1-Hom}_G$ with $\mD\mB\!\on{im_{k,G}}$.
With this modification to Theorem \ref{thm:invariants} and
its proof,\footnote{The last step of the proof is Lemma
  \ref{lem:injective}.} 
we obtain the isomorphism
$$
  \on{HH^\bullet}(A\rtimes G,A\rtimes G) 
  \cong\(\bigoplus_{g\in G}\on{HH^\bullet}(A,Ag)\)^G
$$
of \cite[Prop 3]{Dolgushev:Etingof}.\footnote{It appears that this
  result goes back to an unpublished 
  preprint by Brylinski, dating from 1987.}
\subsection{Coherent sheaves}
Let $X$ be a smooth projective variety over $k$, let $\mC\!\on{oh}(X)$
be the category of coherent sheaves on $X$, and assume that a finite
group $G$ acts on $X$ from the left. Then $g\in
G$ acts on $\mC\!\on{oh}(X)$ via
\begin{equation}
  \label{eq:sheaf-action}
  \mF\longmapsto g_*\mF,
\end{equation}
and the category of equivariant objects in $\mC\!\on{oh}(X)$ is
identified with the category of $G$-equivariant sheaves on $X$,
$$
  \mC\!\on{oh}(X)^G\simeq \mC\!\on{oh_G}(X).
$$
We may view $\mC\!\on{oh_G}(X)$ as the category of sheaves on the orbifold
quotient $[X/G]$. 

It is possible to extend the 2-category of kernels $\mV\!\on{ar}$ of
Section \ref{sec:DCoh-mu} so that it contains both
$X$ and $[X/G]$ as objects 
(see \cite[p.6]{Caldararu:Willerton} or
\cite{Caldararu:III}). 
The formalism remains the same as in Section \ref{sec:DCoh-mu}.
Recall that the $G$-action \eqref{eq:sheaf-action} on
$\mD^+\mC\!\on{oh}(X)$ lifts to a 2-representation
$$
  \longmap{\varrho}{G}{\on{1Hom}(X,X)}
$$
in $\mV\!\on{ar}$, given by the kernels
$
  \varrho(g) = \mathscr O_{\Gamma_g}.
$

Let $Y$ be a second proper, smooth space with $G$-action, and
write $\sigma(s):=\mathscr O_{\Gamma_s}$ for the corresponding
2-representation of kernels on $Y$. Using the projection formula and
flat base-change (c.f.\ \cite[2.7]{Caldararu:derived}), one finds
natural isomorphisms
$$
  \varrho(g)\circ\mathcal K \cong R(1,g)_*\mathcal K
$$
and
$$
  \mathcal K\circ\sigma(s\inv) \cong L(s\inv,1)^*\mathcal K
  \cong R(s,1)_*\mathcal K.
$$
So, 
$$
  \varrho(g)\circ\mathcal K\circ\sigma(s\inv)\cong R(s,g)_*\mathcal K,
$$
and the categories of equivariant 1-homomorphisms of Definition 
\ref{def:1Hom_G} are identified as
$$
  \on{1Hom_G}(Y,X)\simeq \(\mD^+\mC\!\on{oh}(Y\times X)\)^G
$$
and
$$
  \on{1Hom_G^{triv}}(Y,X)\simeq \(\mD^+\mC\!\on{oh}(Y\times X)\)^{G\times G}.
$$
As in the previous section, we need to replace these with 
$$
  \mD^+\mC\!\on{oh_G}(Y\times X)
  \quad\text{and}\quad \mD^+\mC\!\on{oh_{G\times G}}(Y\times X).
$$
With this modification,
$[X/G]$ satisfies the universal property for test-objects $Y$ as above,
$$
  \mD^+\mC\!\on{oh_{G}}(Y\times [X/G])\simeq
  \mD^+\mC\!\on{oh_{G\times G}}(Y\times X),
$$
and the proof of Theorem \ref{thm:invariants} goes through (using
Lemma \ref{lem:injective}) to yield
$$
  \on{HH^\bullet}([X/G]) \cong \bigoplus_{[g]}
  \ttr^\bullet(g)^{C_g}. 
$$
The right-hand side of this isomorphism is identified by Theorem
\ref{thm:ttr(g)}. We have proved
\begin{Thm}
  We have
  $$
    \on{HH^\bullet}([X/G]) \cong \(\bigoplus_{g,\alpha}
    \on{HH^{\bullet-\on{codim}(X^g_\alpha)}}\(X^g_\alpha,\det(N^g)\)\)^{G}, 
  $$
  where $\alpha$ runs over the connected components of the fixed point
  sets. 
\end{Thm}
\subsection{The example of a trivial $G$-action}
%
%
%
%
%
%
Let $\mathcal B$ be a $k$-linear abelian category, and 
assume that $G$
acts trivially (from the left) on $\mathcal B$. Then the data of an
equivariant object $$(B,\epsilon)\in\ob \mathcal B^G$$ are equivalent to
that of an object $B\in \ob\mathcal B$ together with a right
$G$-action 
$$
  \map\epsilon{k[G]}{\on{End}(B)}.
$$
By \cite[5.11]{Deligne}, we have an equivalence of
categories 
\begin{eqnarray*}
  \mathcal B^G&\simeq&\(\on{k[G]-mod^f}\)\widehat\boxtimes\mathcal B\\
               & \simeq & \(W_1\boxtimes \mathcal
               B\)\boxplus\dots\boxplus \(W_r\boxtimes\mB\)\\
              & \simeq &\(\on{k[G]-mod^f}\)\boxtimes\mathcal B,\\
\end{eqnarray*}
where the $W_i$ represent the isomorphism classes of irreducible
representations of $G$. 
Here third equivalence follows from the construction of $\boxtimes$,
and the second equivalence holds, because the category in the second row
is already abelian.
\begin{Exa}
  Let $\mB$ be $\mB=\on{Vect^\CC}(X)$. Then
  $\mB^G=\on{Vect_G^\CC}(X)$ is the category of vector bundles on $X$
  with fibre-preserving $G$-action. 
  In this case the above equivalences boil down to the well
  known fact that any $G$-vector bundle $V$ over a base-space with
  trivial $G$-action can be decomposed as
  $$
    V\cong\bigoplus_i W_i\tensor V_i,
  $$
  where the $V_i\in\on{Vect^\CC}(X)$ are non-equivariant vector bundles.   
\end{Exa}
We have seen two ways to calculate the center of $\mB^G$. By Theorem
\ref{tensor-Thm}, 
$$
  \Cent\(\mB^G\)\cong\on{Center}(k[G])\tensor\Cent(\mB).
$$
On the other hand, since $G$ acts trivially on $\mB$, Theorem
\ref{thm:invariants} becomes
\begin{eqnarray*}
  \Cent\(\mB^G\)& \cong&
  \bigoplus_{[g]}\ttr(g)^{C_g}\\
  &\cong &   \bigoplus_{[g]}\Cent(\mB).
\end{eqnarray*}
Indeed, these two results agree.
In the special case where $\mB=\on{Vect_k}$, we can be more specific
and identify the isomorphism 
$$
  i\negmedspace :
  \on{Center}(k[G])\longrightarrow \bigoplus_{[g]} k 
$$
of Theorem \ref{thm:invariants}: its
inverse $i\inv$ sends the $[g]$th basis vector to the
element 
$$e_{[g]} :=\sum_{h\in C_g}h\inv$$
of $\on{Center}(k[G])$.
\appendix
\section{Clean intersections}
This appendix collects some results about clean intersections that are
probably well-known to the experts. I did not know references, so I am
including proofs.
\begin{Lem}
\label{lem:normal}
  Let $W\sub Z$ be a closed subscheme defined by the ideal sheaf
  $\mI\sub \mathscr O_Z$. Let $w\in W$, and write $\mathfrak m$ for
  the maximal ideal of $\mathscr O_{Z,w}$. Then the conormal space of
  $W$ in $Z$ at $w$ is naturally identified with the quotient
  $$
    N^\vee_w \medspace\cong\medspace \mI_w / \(\mI_w\cap \mathfrak m^2\).
  $$
  If $\mI_w$ is a prime ideal in $\mathscr O_{Z,w}$, this simplifies
  to 
  $$
    N^\vee_w \medspace\cong\medspace \mI_w / \mathfrak m\mI_w.
  $$
\end{Lem}
\begin{Pf}{}
  Write $\mathfrak n$ for the maximal ideal of $$\mathscr
  O_{W,w}=\mathscr O_{Z,w}/\mI_w.$$
  Then $\mathfrak n=\mathfrak m/\mI_w$, and we have the short exact
  sequences
  $$
    \xymatrix{
      0\ar[r]&{\mI_w/\mI_w\cap \mathfrak m^2}\ar[r]\ar[d] &
      \mathfrak m/\mathfrak m^2 \ar[r]\ar[d] &
      \mathfrak n/\mathfrak n^2 \ar[r]\ar[d] & 0\\
       0\ar[r]&{N^\vee_w}\ar[r] &
       T^\vee_{Z,w} \ar[r]&
       T^\vee_{W,w} \ar[r]& 0.
    }
  $$
  Hence the first claim. If $\mI_w$ is prime, then 
  $$
    \mI_w\cap\mathfrak m^2 = \mathfrak m\mI_w,
  $$
  hence the second claim.
\end{Pf}{}
\begin{Cor}
  \label{cor:tangent-isomorphism}
  Let $i\negmedspace : W\hookrightarrow Z$ be the inclusion of an
  irreducible closed subscheme, and let $w\in W$. If
  $$
    \longmap{T_{i,w}}{T_{W,w}}{T_{Z,w}}
  $$
  is an isomorphism then $i$ is an insomorphism in a neighbourhood of
  $w$.
\end{Cor}
\begin{Pf}{}
  The map $T_{i,w}$ is an isomorphism if and only if
  $N^\vee_w=\{0\}$. Since $W$ is irreducible, this is equivalent to
  $I_w = \mathfrak m\mI_w$. By Nakayama's Lemma, this is equivalent to
  $\mI_w=0$. 
\end{Pf}{}
\begin{Lem}
  \label{lem:Y1}
  Let $i\negmedspace : W\hookrightarrow Z$ be the inclusion of a closed
  subscheme, and let $w\in W$. Then there exists an open neighbourhood
  $U$ of $w$ in $Z$ and a regular closed subscheme $Y\sub U$ such that
  $Y$ contains $W\cap U$ and $T_{Y,w}=T_{W,w}$ inside $T_{Z,w}$. 
\end{Lem}
\begin{Pf}{}
  Let $\mathfrak m$, $\mathfrak n$, and $\mI$ be as in Lemma
  \ref{lem:normal}. Pick a basis $f_1,\dots,f_d$ of 
  $$N^\vee_w\medspace \cong\medspace \mI_w/\(\mI_w\cap\mathfrak
  m^2\).$$ 
  By \cite[Cor 4.2.12]{Liu}, this sequence is regular. Hence we can
  find an open neighbourhood $U$ of $w$ in $Z$ and a regular sequence 
  $$
    \tilde f_1,\dots,\tilde f_d \in\Gamma (U,\mI),
  $$
  representing $f_1,\dots,f_d$. By construction, we have $W\sub Y$,
  the closed immersion $Y\hookrightarrow U$ is regular at $w$, and the
  conormal spaces of $W$ and $Y$ inside $Z$ at $w$ agree.
\end{Pf}{}
\begin{Lem}[{\cite[Lem.5.1]{LiLi:wonderful}}]
  Let $Z$ be a nonsingular variety, and let $X$ and $Y$ be nonsingular
  closed subvarieties of $Z$. Then the connected components of the 
  scheme-theoretic intersection
  $X\cap Y$ of $X$ and $Y$ in $Z$ are non-singular if and only if the
  following differential geometric conditions are satisfied:
  \begin{enumerate}
  \item the connected components of the ``set-theoretic intersection''
    $$(X\cap Y)_{red}$$
    are nonsingular varieties, and
  \item we have $$T_X\cap T_Y=T_{X\cap Y}$$ inside $T_Z$.
  \end{enumerate}  
\end{Lem}
In the context of differential geometry, these conditions define the
notion of {\em clean intersection}, due to Bott \cite{Bott:56}.
\begin{Lem}\label{lem:clean}
  In the situation of the previous Lemma the canonical map of vector
  bundles over $X\cap Y$ 
  $$
    T_X+T_Y\hookrightarrow T_{X\cup Y}
  $$
  is an isomorphism. Here the sum on the left-hand side is taken
  inside $T_Z{}\vert_{X\cap Y}$.
\end{Lem}
\begin{Pf}{}
  Fix $z\in X\cap Y$, let $\mathfrak m\subset\mathscr O_{Z,z}$ be the
  maximal ideal, and let $\mI$ and $\mathcal J$ denote the stalks of the
  ideal sheaves of $X$ and $Y$ at $z$. Then we have a cartesian square
  $$
    \xymatrix{
    \mI\cap \mathcal J/\mI\cap\mathcal J\cap\mathfrak m^2 \ar@{>->}[r]\ar@{>->}[d]  &
    \mI/\mI\cap\mathfrak m^2 \ar@{>->}[d]  \\
    \mathcal J/\mathcal J\cap\mathfrak m^2 \ar@{>->}[r]  &
    \mI+ \mathcal J/(\mI+\mathcal J)\cap\mathfrak m^2. 
    }
  $$
  Lemma \ref{lem:normal} identifies the vector spaces in this diagram
  with the conormal spaces of $X\cup Y$, $X$, $Y$, and $X\cap Y$
  inside $Z$ at $z$. A dimension count now proves the claim.
\end{Pf}{}
\begin{Lem}\label{lem:Y2}
  Let $Z$ be a nonsingular variety, and let $X$, $Y$, and $W$ be
  nonsingular closed subvarieties of $Z$ such that $W$ is contained in
  $X$ and in $Y$ (``contained'' means as a closed subscheme). Let
  $w\in W$, and assume that the map
  $$
    T_{W,w} \hookrightarrow T_{X,w}\cap T_{Y,w}
  $$
  is an isomorphisms. (Here the intersection on the right-hand side is
  inside $T_{Z,z}$.) Then there exists a neighbourhood of $w$ in $Z$
  inside which the map $i\negmedspace : W\hookrightarrow X\cap Y$ is an
  isomorphism. 
\end{Lem}
\begin{Pf}{}
  We know that the composite
  $$
    \xymatrix{
    T_{W,w}\ar@{>->}[0,2]^{T_{i,w}}&&T_{X\cap
      Y}\ar@{>->}[0,2]&&T_{X}\cap T_Y 
    }
  $$
  is an isomorphism. Hence $T_{i,w}$ is an isomorphism, as
  well. Applying Corollary \ref{cor:tangent-isomorphism}, we obtain
  the claim.
\end{Pf}{}
We will also need a lemma from homological algebra.
\begin{Lem}\label{lem:injective}
  Let $\mA$ be an abelian category, and assume that $\mathcal A^G$ has
  enough injectives. Let $\mF$ and $\mG$ be
  equivariant objects in $\mA$. Then we have
  $$
    \Hom^\bullet_{\mathcal D^+_G(\mA)}(\mathcal F,\mathcal G)\cong 
    \Hom^\bullet_{\mathcal D^+(\mA)}(\mathcal F,\mathcal G)^G.
  $$  
\end{Lem}
\begin{Pf}{}
  Choose an injective resolution $\mG\to\mathcal I^\bullet$ in
  $\mA^G$. Since the forgetful functor $\map U{\mA^G}\mA$ is a
  right-adjoint, it preserves injectives. Therefore, we have
  \begin{eqnarray*}
    \Hom^\bullet_{\mathcal D^+_G(\mA)}(\mathcal F,\mathcal G) &\cong
    &
       \Hom^\bullet_{\mathcal Ch^+_G(\mA)}(\mathcal F,\mathcal
       I^\bullet)\\
   &\cong&  
       \Hom^\bullet_{\mathcal Ch^+(\mA)}\(\mathcal F,\mathcal
       I^\bullet\)^G\\
   &\cong&  
       \Hom^\bullet_{\mathcal D^+(\mA)}\(\mathcal F,\mathcal
       G\)^G.
  \end{eqnarray*}
\end{Pf}{}
\bibliographystyle{alpha}
\bibliography{innerproducts}
\end{document}